\documentclass[aap,citesort,MSNbibl,dvips]{arximspdf}
\usepackage{graphics}
\doi{10.1214/08-AAP575}
\volume{19}
\issue{4}
\pubyear{2009}
\firstpage{1292}
\lastpage{1318}

\makeatletter

\newtheorem{theorem}{Theorem}[section]
\newtheorem{lemma}[theorem]{Lemma}
\newtheorem{prop}[theorem]{Proposition}

\newproclaim{remarks}{Remarks}
\newproclaim{remark}{Remark}

\makeatother

\begin{document}
\begin{frontmatter}

\title{Subcritical regimes in some models of continuum~percolation}
\runtitle{Subcritical regimes in continuum percolation}

\begin{aug}
\author[A]{\fnms{Jean-Baptiste}
\snm{Gou\'{e}r\'{e}}\corref{}\ead[label=e1]{jbgouere@univ-orleans.fr}}
\runauthor{J.-B. Gou\'{e}r\'{e}}
\affiliation{Universit\'{e} d'Orl\'{e}ans}
\address[A]{Universit\'{e} d'Orl\'{e}ans\\
MAPMO\\
B.P. 6759\\
45067 Orl\'{e}ans Cedex 2\\
France\\
\printead{e1}} 
\end{aug}
\pdfauthor{Jean-Baptiste Gouere}

\received{\smonth{4} \syear{2008}}
\revised{\smonth{10} \syear{2008}}

%
\begin{abstract}
We consider some continuum percolation models. We are mainly interested
in giving some sufficient conditions for the absence of percolation. We
give some general conditions and then focus on two examples. The first
one is a multiscale percolation model based on the Boolean model. It
was introduced by Meester and Roy and subsequently studied by
Menshikov, Popov and Vachkovskaia. The second one is based on the
stable marriage of Poisson and Lebesgue introduced by Hoffman, Holroyd
and Peres and whose percolation properties have been studied by Freire,
Popov and Vachkovskaia.
\end{abstract}

%

\end{frontmatter}

\section{Introduction and statement of the main results}

\subsection{Introduction}
\label{introintro}

In this paper, we study some continuum percolation models.
We are mainly interested in giving some sufficient conditions for the
absence of percolation.
We give some general conditions and then apply them to two examples:
multiscale percolation and the stable marriage of Poisson and Lebesgue.
The aim of Section \ref{introintro} is to give a quick description of
these two examples.
We begin by recalling the Boolean model.

\subsubsection*{The Boolean model}
We center a ball of random radius at each point of a homogeneous
Poisson point process of the Euclidean space $\mathbb{R}^d$, $d\ge2$.
We assume that the radii of the balls are independent copies of a given
positive random variable $R$.
We also assume that the radii are independent of the point process.
We denote by $\lambda$ the density of the Poisson point process.
We denote by $\Sigma(\lambda)$ the union of the balls,
by $S(\lambda)$ the connected component of $\Sigma(\lambda)$ that
contains the origin
and by $D(\lambda)$ the Euclidean diameter of $S(\lambda)$.

When $R$ is bounded, there exists a sharp phase transition
(see, e.g., \cite{MeesterRoylivre}, Section~12.10 in \cite
{Grimmettpercolation} when $R=1$
or the papers \cite{MeesterRoySarkar,ZuevI,MenshikovSidorenkocoincidence} and \cite{ZuevII}):
if the density~$\lambda$ of the point process is below a critical
value $\lambda_c>0$,
then $S(\lambda)$ is almost surely bounded and its diameter $D(\lambda
)$ admits exponential moments;
whereas if $\lambda$ is above~$\lambda_c$ then $S(\lambda)$ is
unbounded with positive probability.

The case where $R$ is unbounded was studied by Hall in \cite
{Hallcontinuumpercolation}
(see also the books \cite{Halllivre} and \cite{MeesterRoylivre}).
Hall proved that if $E(R^{2d-1})$ is finite,
then the set $S(\lambda)$ is almost surely bounded for small enough
$\lambda$.
If $E(R^d)$ is infinite, then such behavior does not happen:
whatever the value of the density $\lambda$, the set $\Sigma(\lambda
)$ is almost surely the whole space.
The latter result still holds when the underlying point process is only
assumed to be an almost surely nonempty and stationary point process
(see \cite{MeesterRoylivre}, Proposition 7.3).
In \cite{Gbooleanmodel}, we proved that the set $S(\lambda)$ is
almost surely bounded for small enough $\lambda$ if and only if
$E(R^d)$ is finite.
We also proved that, for any $s>0$, $E(D(\lambda)^s)$ is finite for
small enough $\lambda$ if and only if $E(R^{d+s})$ is finite.
We refer to \cite{Gbooleanmodel} for further bibliographical information.
The idea developed in \cite{Gbooleanmodel} for the Boolean model can
be developed further to investigate the following models.

\subsubsection*{A multiscale percolation model}
We refer to Section \ref{s-multi} for details.
We keep the objects defined in the previous paragraph.
Let ${(\Sigma_n(\lambda))}_{n\ge0}$ be a sequence of
independent\vspace*{1pt}
copies of $\Sigma(\lambda)$.
Let $a>1$ be a scale factor.
We define a new random set $\widetilde{\Sigma}(\lambda,a)$ by
\[
\widetilde{\Sigma}(\lambda,a)=\bigcup_n a^{-n} \Sigma_n(\lambda).
\]
We are interested in properties of the connected components of
$\widetilde{\Sigma}(\lambda,a)$.
We say that the model is subcritical if the connected components can be
small (see Section~\ref{s-multi} for precise statements).

This model was introduced by Meester and Roy in \cite{MeesterRoylivre}.
They considered the case where the radius $R=1$ and the dimension $d=2$.
They proved that $\widetilde{\Sigma}(\lambda,a)$ is in a subcritical
phase as soon as one of the following conditions holds:
\begin{enumerate}
\item the density $\lambda$ is small enough;
\item the density $\lambda$ is such that $\Sigma(\lambda)$ is in the
subcritical phase and $a$ is large enough.
\end{enumerate}
%

In \cite{Menshikovalmulti}, Menshikov, Popov and Vachkovskaia
considered the case where the dimension $d$ is arbitrary.
They proved that if $\Sigma(\lambda)$ is in the subcritical phase
then, for large enough $a$, $\widetilde{\Sigma}(\lambda,a)$ is also in
a subcritical phase.

In \cite{Menshikovalmultiunbounded}, Menshikov, Popov and
Vachkovskaia studied the case where $R$ is random.
Let us emphasize that they did not assume $R$ to be bounded.
They considered the following condition:
%
%
\begin{equation} \label{c-decay}
P\bigl(D(\lambda) \ge r\bigr)r^d \to0 \qquad\mbox{as } r\to\infty.
\end{equation}
[Let us recall that $D(\lambda)$ denotes the Euclidean diameter of the
connected component of $\Sigma(\lambda)$ containing $0$.]
Under some further technical conditions,
they proved that, if~(\ref{c-decay}) holds, then for large enough $a$,
$\widetilde{\Sigma}(\lambda,a)$ is in a subcritical phase.
This is a generalization of the previous result.
Indeed, when $R$ is bounded, (\ref{c-decay}) holds as soon as $\Sigma
(\lambda)$ is in the subcritical phase.
When $R$ is unbounded, one can make the following remarks about the
conditions under which there exists $\lambda>0$ such that~(\ref
{c-decay}) holds.
Let $\varepsilon>0$.
Condition (\ref{c-decay}) holds as soon as $E(D(\lambda
)^{d+\varepsilon
})$ is finite.
Therefore, by the result of \cite{Gbooleanmodel} previously cited,
there exists $\lambda>0$ such that (\ref{c-decay}) holds as soon as
$E(R^{2d+\varepsilon})$ is finite.
On the other hand, if (\ref{c-decay}) holds then $E(D(\lambda
)^{d-\varepsilon})$ is finite and therefore $E(R^{2d-\varepsilon})$
is finite.
Therefore, the existence of $\lambda>0$ such that (\ref{c-decay})
holds is roughly equivalent to the finiteness of $E(R^{2d})$.

In this paper, we prove the following result in which $a>1$ is fixed:
$\widetilde{\Sigma}(\lambda,a)$ is in a subcritical phase for small
enough $\lambda$ if and only if $E(R^d\max(\ln(R),0))$ is finite.
This is a corollary of Theorem \ref{th-poisson-non-percolation}, which
is one of our main abstract results.

\subsubsection*{Stable marriage of Poisson and Lebesgue}
We refer to Section \ref{s-marriage} for details.
The following model was introduced by Hoffman, Holroyd and Peres in
\cite{HoffmanHolroydPeresmarriage}.
Let $\alpha>0$ be a parameter called appetite.
Let $\chi$ be a homogeneous Poisson point process with density $1$ on
$\mathbb{R}^d$.
In \cite{HoffmanHolroydPeresmarriage},
the authors showed that there was essentially a unique way to give in a
stable way to points of $\chi$
disjoint territories of $\mathbb{R}^d$ of volume at most $\alpha$.
We defer the definition of stability to Section \ref{s-marriage}.
Very roughly, it means that the distances between points of $\chi$ and
points of their territories are minimal.



In this paper, we are interested in percolation properties of the union
$T(\alpha)$ of all the territories.
Let $S(\alpha)$ denote the connected component of $T(\alpha)$ that
contains the origin.
In \cite{FreirePopovVachkovskiapercolationmarriage}, Freire, Popov
and Vachkovskaia proved, among other things, that~$S(\alpha)$ was
almost surely bounded for small enough $\lambda$.
In this paper, we prove the following stronger result,
in which $D(\alpha)$ denotes the Euclidean diameter of~$S(\alpha)$.
For small enough $\lambda$, for all $n\ge0$, $E(D(\alpha)^n)$ is finite.

To prove this result, we first show that $T(\alpha)$ is dominated by a
dependent percolation process.
This was already the first step in the proof of \cite
{FreirePopovVachkovskiapercolationmarriage}.
We then apply to this dependent percolation process Theorem \ref
{th-cs-non-percolation},
which is the main abstract theorem of our paper.

\subsection{Some notation}
\label{s-notations}

For the whole of the paper, we fix an integer $d\ge1$.
Let~$|\cdot|$ be the Lebesgue measure on $\mathbb{R}^d$.
We denote by $\|\cdot\|$ the Euclidean norm on $\mathbb{R}^d$,
by $B(x,r)$ the open Euclidean ball centered at $x\in\mathbb{R}^d$ with
radius $r\ge0$
and by $\overline{B}(x,r)$ the closed Euclidean ball centered at $x\in
\mathbb{R}^d$ with radius $r\ge0$.

When a point process $\xi$ on $\mathbb{R}^d\times\,]0,+\infty[$ is
given we
define the following objects.
We let
\[
\Sigma=\bigcup_{(c,r)\in\xi} B(c,r).
\]
(When we write $(c,r)\in\xi$ we implicitly assume that $c$ belongs to
$\mathbb{R}^d$ and that~$r$ belongs to $]0,+\infty[$.)
We denote by $S$ the connected component of $\Sigma$ which contains $0$.
(We let $S=\varnothing$ if $0$ does not belong to $\Sigma$.)
We define a random variable $M$ as follows:
%
%
\begin{equation} \label{defM}
M=\sup_{x\in S} \|x\|.
\end{equation}
(We let $M=0$ if $S$ is empty.)
We say that percolation occurs if $S$ is unbounded:
\[
\{\mbox{percolation}\}=\{S\mbox{ is unbounded}\}.
\]

\subsection{Boolean model induced by Poisson point processes}
\label{s-poisson}

Let $\lambda>0$ and let~$\mu$ be a locally finite measure on
$]0,+\infty[$.
Let $\xi$ be a Poisson point process on $\mathbb{R}^d\times
\,]0,+\infty[$
whose intensity measure is the product of $\lambda|\cdot|$ and $\mu$.
We denote by $P_{\lambda,\mu}$, $E_{\lambda,\mu}$ the associated
probability measure and expectation, respectively.
As distinct points of $\xi$ have distinct coordinates on $\mathbb
{R}^d$, we
can write
\[
\xi=\{(c,r(c)), c\in\chi\},
\]
where $\chi$ denotes the projection of $\xi$ on $\mathbb{R}^d$.
If the measure $\mu$ is a probability measure then
$\chi$ is a Poisson point process on $\mathbb{R}^d$ whose intensity is
$\lambda|\cdot|$.
Moreover, under this assumption, if we condition on $\chi$ then the
$r(c), c\in\xi$ are i.i.d. with common distribution $\mu$.
(We shall not use this result.)
We refer to \cite{Kallenbergrandommeasures,Moller,Neveupp} for
background on point processes
and to \cite{Halllivre,MeesterRoylivre} for Boolean models.

We prove the following results:
\begin{theorem} \label{th-poisson-non-percolation}
Assume $d \ge2$.
There exists $\lambda_0>0$ such that $P_{\lambda,\mu}(\mbox{perco}\-\mbox{lation})=0$ for all $\lambda\in\,]0,\lambda_0[$
if and only if the following assertions hold:
\begin{longlist}[A2.]
\item[A1.] The supremum ${\sup_{r>0} r^d
\mu([r,+\infty[)}$ is finite.
\item[A2.] The integral ${\int_{[1,+\infty
[} \beta^d \mu(d\beta)}$ is finite.
\end{longlist}

If $d=1$, then assumptions \textup{A1} and \textup{A2}
together are sufficient conditions;
assumption \textup{A2} is a necessary condition.
\end{theorem}

\begin{remarks*}
\begin{enumerate}
\item For all $\rho>1$, assumption A1 is equivalent to
the following one:
\begin{eqnarray*}
\sup_{r>0} r^d\mu([r,\rho r])<\infty
\end{eqnarray*}
(see Lemma \ref{l-Abis}).
Notice that the probability of $0$ belonging to a ball of the process
with radius in $[r,\rho r]$ is
\[
P_{\lambda,\mu} \biggl(0 \in\bigcup_{(c,\beta)\in\xi \dvtx  \beta
\in[r,\rho r]} B(c,\beta) \biggr) = 1-\exp \biggl(-\lambda\int
_{[r,\rho r]} \mu(d\beta) |B(0,\beta)| \biggr).
\]
Assumption A1 therefore means that those probabilities
are bounded away from~$1$.
%
\item Assume in this remark that $\mu$ is a finite measure.
Then, the integral $\int_{]0,1[} \beta^d\times \mu(d\beta)$ is finite.
Therefore, assumption A2 holds if and only if the integral
%
%
\begin{equation} \label{courseok}
\int_{]0,+\infty[} \beta^d \mu(d\beta)
\end{equation}
is finite.
Moreover, as assumption A1 holds as soon as the integral
(\ref{courseok}) is finite,
the assumptions A1 and A2 together are
also equivalent to the finiteness of the integral (\ref{courseok}).
\end{enumerate}
\end{remarks*}
\begin{theorem} \label{th-poisson-moment-fini}
Let $s>0$ be a positive real.
Assume $d\ge2$.
There exists $\lambda_0>0$ such that $E_{\lambda,\mu}(M^s)$ is
finite for all $\lambda\in\,]0,\lambda_0[$
if and only if the following assertions hold:
\begin{longlist}[A3.]
\item[A1.] The supremum ${\sup_{r>0} r^d
\mu([r,+\infty[)}$ is finite.
\item[A3.] The integral ${\int_{[1,+\infty
[} \beta^{d+s} \mu(d\beta)}$ is finite.
\end{longlist}

If $d=1$, then assumptions \textup{A1} and \textup{A3}
together are sufficient conditions;
assumption \textup{A3} is a necessary condition.
\end{theorem}
\begin{remark*}
If $\mu$ is a finite measure, then assumptions A1
and A3 together are equivalent to the finiteness
of the integral $\int_{]0,+\infty[} \beta^{d+s} \mu(d\beta)$.
\end{remark*}

Theorems \ref{th-poisson-non-percolation} and \ref
{th-poisson-moment-fini} are essentially consequences of Theorem \ref
{th-cs-non-percolation} stated in the next subsection.
Theorems \ref{th-poisson-non-percolation} and \ref
{th-poisson-moment-fini} are generalizations of the main
results of~\cite{Gbooleanmodel} in which $\mu$ is assumed to be a
finite measure.

The proofs are given in Section \ref{s-preuve-poisson}.

\subsection{Boolean model induced by more general point processes}

\label{s-general}

Let $\xi$ be a point process on $\mathbb{R}^d\times\,]0,+\infty[$.
We assume that the law of $\xi$ is invariant under the action of the
translations of $\mathbb{R}^d$:
for all $t\in\mathbb{R}^d$, the point processes $\{x-(t,0), x\in\xi
\}$ and
$\xi$ have the same law.
We also assume that the intensity measure of $\xi$ is locally finite.
Therefore, the intensity measure of $\xi$ is the product of the
Lebesgue measure on~$\mathbb{R}^d$
by a locally finite measure on $]0,+\infty[$ that we denote by $\mu$.

The main result of this paper is the following theorem.
\begin{theorem} \label{th-cs-non-percolation}
Let $C>0$.
There exists $D>0$, that depends only on $d$ and~$C$, such that the
following hold.

If the following properties are fulfilled:
\begin{longlist}[B0.]
\item[B0.] for all $r>0$ and all $x\in\mathbb{R}^d\setminus B(0,C r)$
the point processes
\[
\xi\cap B(0,r) \times\,]0,r] \quad\mbox{and}\quad \xi\cap B(x,r) \times\,]0,r]
\]
\end{longlist}\vspace*{-3pt}
are independent;\vspace*{-3pt}
\begin{longlist}[B1.]
\item[B1.] ${\sup_{r>0} r^d \mu([r,+\infty[)}
\le D$;\vspace*{0.5pt}
\item[B2.] the integral ${\int_{[1,+\infty[}
\beta^d \mu(d\beta)}$ is finite,\vspace*{-3pt}
\end{longlist}
then the set $S$ is almost surely bounded.
Let $s$ be a positive real.
If, moreover,\vspace*{-3pt}
\begin{longlist}[B3.]
\item[B3.] ${\int_{[1,+\infty[} \beta^{d+s}
\mu(d\beta)<\infty}$,\vspace*{-3pt}
\end{longlist}
then $E(M^s)$ is finite.
\end{theorem}

\begin{remarks*}
\begin{enumerate}
\item
The independence assumption B0 is fulfilled if $\xi$ is
a Poisson point process and $C \ge2$.
\item We give a strenghtened version of Theorem \ref
{th-cs-non-percolation} in Section \ref{s.existence}
(see Theorems~\ref{th-cs-non-percolation-I}, \ref
{th-cs-non-percolation-II} and \ref{th-cs-non-percolation-III}).
In those theorems, the independence assumption is weakened and the
conclusions are strenghtened.
\end{enumerate}

The proof is given in Section \ref{s-preuve-generale}.
We give some ideas of the proof in Section \ref{s-key} after the
statement of key Proposition \ref{laprop}.
\end{remarks*}

\subsection{Multiscale percolation model}
\label{s-multi}

Let $\lambda>0$ and $\nu$ be a probability measure on $]0,+\infty[$.
We make the following assumption:
%
%
\begin{equation} \label{hyp-moment-fractal}
\int_{]0,+\infty[} r^d \nu(dr)<\infty.
\end{equation}
Let ${(\xi_n)}_{n\ge0}$ be a sequence of independent Poisson point
processes on $\mathbb{R}^d\times\,]0,\break+\infty[$ whose intensity
is the product of $\lambda|\cdot|$ by $\nu$.
Let $a>1$.
We define a new point process by:
\[
\xi= \bigcup_{n\ge0} a^{-n}\xi_n.
\]

\begin{lemma} \label{l-fractal}
The point process $\xi$ is a Poisson point process whose intensity is
the product of $\lambda|\cdot|$ by the locally finite measure $\mu$
on $]0,+\infty[$ defined by:
%
%
\begin{equation} \label{f-def-mu}
\mu(B)=\sum_{n \ge0} a^{nd} \nu(a^n B).
\end{equation}
\end{lemma}

As in Section \ref{s-notations} we associate with $\xi$ two random
sets $\Sigma$ and $S$.
We denote by~$P_{\lambda,\nu}^a$ the associated probability measure.
We also denote by $\Sigma_n$ the random sets associated with the
processes $a^{-n}\xi_n$.

\begin{remarks*}
\begin{enumerate}
\item
For all integer $n\ge1$, $\Sigma_n$ is an independent copy of
$a^{-n}\Sigma_0$.
\item If (\ref{hyp-moment-fractal}) is not fulfilled then, for all
$\lambda>0$,
percolation occurs with positive probability in $\Sigma_0$ (by Theorem
\ref{th-poisson-non-percolation}) and then in $\Sigma$.
Actually, by Lemma \ref{l-perco-grand}, if (\ref{hyp-moment-fractal})
is not fulfilled then, for all $\lambda>0$, $\Sigma_0=\mathbb{R}^d$
almost surely.
Therefore, assumption~(\ref{hyp-moment-fractal}) is not a restriction.
\item One can easily check that $0$ belongs almost surely to $\Sigma$.
Therefore, the Lebesgue measure of the complement of $\Sigma$ is
almost surely $0$.
We will nevertheless see that the connected components of $\Sigma$ can
be bounded.
\end{enumerate}

This model was introduced by Meester and Roy in a two-dimensional
setting in~\cite{MeesterRoylivre}.
Let us denote by $\delta_1$ the Dirac mass at $1$.
Let us say that the event $\{$left--right crossing$\}$ occurs if
$[0,1]^2\setminus\Sigma$ contains a connected component which
intersects left- and right-hand sides of $[0,1]^2$.
Let us denote by $\lambda_c$ the critical density for the Boolean
model when all radii equal $1$.
(Thus, if $\lambda<\lambda_c$, the connected components of the
$\Sigma_n$ are almost surely bounded;
whereas if $\lambda>\lambda_c$, this is not the case.)
Meester and Roy proved the following result, in which the radii of the
unscaled process $\Sigma_0$ equal $1$.
\end{remarks*}
\begin{theorem}[(\cite{MeesterRoylivre})]
Assume $d=2$.
\begin{enumerate}
\item Let $a>1$. If $\lambda>0$ is small enough, then $P_{\lambda
,\delta_1}^a(\mbox{left--right crossing})$ is positive.
\item Let $\lambda<\lambda_c$. If $a$ is large enough, then
$P_{\lambda,\delta_1}^a(\mbox{left--right crossing})$ is positive.
\end{enumerate}
\end{theorem}

In \cite{Menshikovalmulti}, Menshikov, Popov and Vachkovskaia
considered the case where the dimension $d$ is arbitrary
and the radii of the unscaled process $\Sigma_0$ equal $1$.
They proved the following result.
\begin{theorem}[(\cite{Menshikovalmulti})]\label{th-MPV1}
Assume $d\ge2$.
If $\lambda<\lambda_c$ then, for all $a$ large enough,
\[
P_{\lambda,\delta_1}^a(S \mbox{ is bounded})=1.
\]
\end{theorem}

The ideas of their proof are the following.
(Those ideas are used in their paper through a discretization of space;
we describe them in a slightly more geometric way.)
Assume that $C$ is a connected component of $\Sigma_n \cup\Sigma
_{n+1}$ whose diameter is at least $\alpha a^{-n}$ for a given $\alpha>0$.
Then, $C$ is included in the union of the following kind of sets:
\begin{enumerate}
\item connected components of $\Sigma_{n+1}$ whose diameter is at
least $\alpha a^{-n}$;
\item balls of $\Sigma_n$ enlarged by the factor $1+\alpha$ [same
centers but the radii are $(1+\alpha)a^{-n}$ instead of $a^{-n}$].
\end{enumerate}
Then, they show that the union of all those sets is stochastically
dominated by a Boolean model
where all radii equals $(1+\alpha)a^{-n}$ ($1+\alpha$ times those of
$\Sigma_n$)
and where the density of the set of centers is $(1+\alpha')a^{nd}$ for
a suitable $\alpha'>0$ ($1+\alpha'$ times the corresponding density
for $\Sigma_n$).
The proof of this fact relies partly on the exponential decay of the
size of the components in the subcritical phase.
In some sense, one can therefore control percolation in the union of
two models by percolation in one model.
Iterating the argument with some care in the constants~$\alpha$ and
$\alpha'$,
one sees that one can control percolation in the multiscale model by
percolation in a subcritical model.
This yields the result.

In \cite{Menshikovalmultiunbounded} the same authors considered the
case where the radii are random and unbounded.
Let us define $\Theta$ by
\[
\Theta=\{\lambda>0 \dvtx  P_{\lambda,\nu}(D_0>n)n^d \to0 \mbox{ as }
n\to\infty\},
\]
where $D_0$ denotes here the diameter of the connected component of
$\Sigma_0$ containing~$0$.
Let $\widetilde{\lambda}_c$ denote the supremum of $\Theta$.
They proved the following generalization of Theorem \ref{th-MPV1}.
\begin{theorem}[(\cite{Menshikovalmultiunbounded})]\label{ooo}
Assume $d\ge2$ and the following:
\begin{enumerate}
\item the set $\Theta$ is nonempty (and thus $\widetilde{\lambda}_c$
is positive);
\item the measure $\nu$ satisfies
\[
\lim_{a\to\infty} \sup_{r\ge1/2} \frac{a^d\nu([ar,+\infty
[)}{\nu([r,+\infty[)}=0
\]
with the convention $0/0=0$.
\end{enumerate}
Then, for all $\lambda< \widetilde{\lambda}_c$, for all large enough
$a$, $P_{\lambda,\nu}^a(S \mbox{ is bounded})=1$.
\end{theorem}
\begin{remarks*}
\begin{enumerate}
\item
When $R=1$, $P_{\lambda,\delta_1}(D_0>n)$ decays exponentially
as soon as $\lambda<\lambda_c$.
Therefore $\widetilde{\lambda}_c=\lambda_c$, in that case.
Theorem \ref{ooo} is thus a generalization of Theorem \ref{th-MPV1}.
\item The assumption $\lambda< \widetilde{\lambda}_c$ is used where,
in the proof of Theorem \ref{th-MPV1},
the exponential decay of the size of the connected components in the
subcritical phase were used.
\item As explained in the \hyperref[introintro]{Introduction} (see Section \ref{introintro}),
the first assumption means roughly that the integral $\int r^{2d} \nu
(dr)$ is finite.
\end{enumerate}

By Theorem \ref{th-poisson-non-percolation}, we easily get the
following result.
\end{remarks*}
\begin{theorem} \label{th-fractal}
There exists $\lambda_0>0$ such that $P_{\lambda,\nu}^a(S \mbox{ is
bounded})=1$ for all $\lambda\in\,]0,\lambda_0[$
if and only if the integral
\[
\int_{[1,+\infty[} \beta^d \ln(\beta) \nu(d\beta)
\]
is finite.
\end{theorem}

The proof is given in Section \ref{s-preuve-fractal}.

\begin{remark*} We can get a similar result about the finiteness of moments of
the diameter of $S$ by Theorem \ref{th-poisson-moment-fini}.
\end{remark*}


\subsection{Stable marriage of Poisson and Lebesgue}
\label{s-marriage}

The following model was introduced in \cite
{HoffmanHolroydPeresmarriage} by Hoffman, Holroyd and Peres.
Let $\chi$ be a locally finite subset of~$\mathbb{R}^d$.
We call the elements of $\mathbb{R}^d$ sites and the elements of $\chi
$ centers.
Let $\alpha\in\,]0,\infty[$ be a parameter called the appetite.
An allocation of $\mathbb{R}^d$ to $\chi$ with appetite $\alpha$ is a
measurable function
\[
\psi\dvtx \mathbb{R}^d \to\chi\cup\{\infty,\Delta\}
\]
such that $|\psi^{-1}(\Delta)|=0$,
and $|\psi^{-1}(a)|\leq\alpha$ for all $a\in\chi$.
We call $\psi^{-1}(a)$ the territory of the center $a$.
A center $a\in\chi$ is sated if $|\psi^{-1}(a)|=\alpha$ and unsated
otherwise.
A site $x\in\mathbb{R}^d$ is claimed if $\psi(x)\in\chi$ and
unclaimed if
$\psi(x)=\infty$.
The allocation is undefined at $x$ if $\psi(x)=\Delta$.

The following definition, given in \cite
{HoffmanHolroydPeresmarriage}, is an adaptation of that introduced by
Gale and Shapley \cite{GaleShapleymarriage}.
Let $a$ be a center and let $x$ be a site with $\psi(x)\notin\{
a,\Delta\}$.
We say that $x$ desires $a$ if
\[
\|x-a\| < \|x-\psi(x)\| \quad\mbox{or}\quad \mbox{$x$ is unclaimed.}
\]
We say that $a$ covets $x$ if
\[
\|x-a\| < \|x'-a\| \qquad\mbox{for some $x'\in\psi^{-1}(a)$,}\quad \mbox{or}\quad \mbox{$a$ is unsated.}
\]
We say that a site-center pair $(x,a)$ is unstable for the allocation
$\psi$ if $x$ desires $a$ and $a$ covets $x$.
An allocation is stable if there are no unstable pairs.

We now assume that $\chi$ is a translation invariant Poisson point
process on $\mathbb{R}^d$.
We assume that its intensity measure is the Lebesgue measure.
(We can see by scaling arguments that there is no loss of generality in
this assumption.)
In \cite{HoffmanHolroydPeresmarriage} it was proved, among other
things, that for any such process
there exists a.s. a $|\cdot|$-a.e. unique stable allocation $\psi$
from $\mathbb{R}^d$ to $\chi$.
Furthermore we have the following phase transition phenomenon:
\begin{enumerate}
\item If $\alpha<1$ (subcritical) then a.s. all centers are sated but
there is an infinite volume of unclaimed sites.
\item If $\alpha=1$ (critical) then a.s. all centers are sated and
$|\cdot|$-a.a. sites are claimed.
\item If $\alpha>1$ (supercritical) then a.s. not all centers are
sated but $|\cdot|$-a.a. sites are claimed.
\end{enumerate}

Let $\mathcal{C}$ be the closure of the union of all territories
\[
\mathcal{C}=\overline{\psi^{-1}(\chi)}.
\]
In \cite{FreirePopovVachkovskiapercolationmarriage}, Freire, Popov
and Vachkovskaia proved, among other things, the following result:
\begin{theorem}[(\cite{FreirePopovVachkovskiapercolationmarriage})]
\label{th-marriage-non-perco}
If $\alpha$ is small enough, then a.s. there is no percolation
in~$\mathcal{C}$.
\end{theorem}

Let $D$ be the diameter of the connected component of $\mathcal{C}$
that contains the origin.
In this paper we give the following consequence of Theorem \ref
{th-cs-non-percolation}:
\begin{theorem} \label{th-marriage-non-perco-moments}
If $\alpha$ is small enough, then for all $s>0$, $E(D^s)$ is finite.
\end{theorem}
\begin{remark*} We must admit that we have not checked out the measurability of $D$.
Actually, we prove that, for small enough $\alpha$, $D^s$ is bounded
above by an integrable random variable.
\end{remark*}

In order to prove Theorem \ref{th-marriage-non-perco}, we first define
a process that dominates the previous one.
This relies on an idea that appeared in \cite
{HoffmanHolroydPeresmarriagetail} [see the proof of Proposition~11(ii)]
and that is used with the same purpose as ours in \cite
{FreirePopovVachkovskiapercolationmarriage} (see Lemma 2.1).
For all $a\in\chi$ we define $R(a,\chi)$ by
\[
R(a,\chi)=\inf\bigl\{r \ge0 \dvtx  \alpha\operatorname{card}\bigl(\chi\cap
\overline
{B}(a,2r)\bigr) \le|B(a,r)|\bigr\}.
\]
We let $R(a,\chi)=\infty$ if there is no such $r$.
We assume henceforth that $\alpha$ is strictly smaller than $2^{-d}$.
This ensures that, almost surely, all the $R(a,\chi)$ are finite
(see Lemma \ref{l-marriage-queue-mu} for a stronger statement).
We can also check that all the $R(a,\chi)$ are positive.
We then define a point process $\xi$ on $\mathbb{R}^d\times
\,]0,+\infty[$ by
\[
\xi=\{(a,2R(a,\chi)), a\in\chi\}.
\]
As in Section \ref{s-notations}, we associate with this process a
random set $\Sigma$.
We have:
\begin{lemma} \label{l-marriage-domination}
For all $\alpha\in\,]0,2^{-d}[$, the set $\mathcal{C}$ is almost
surely contained in the set $\Sigma$.
\end{lemma}

It is therefore sufficient to study the percolation properties of
$\Sigma$.
Theorem \ref{th-marriage-non-perco-moments} follows from an
application of Theorem \ref{th-cs-non-percolation} to the process $\xi$.
A full proof is given in Section \ref{s-proof-marriage}.


\section[Proof of Theorem 1.3 (models induced by general
processes)]{Proof of Theorem \protect\ref{th-cs-non-percolation}
(models induced by general processes)}
\label{s-preuve-generale}

\subsection{Some further notation}

For the whole of the section,
we fix a point process~$\xi$.
We assume that $\xi$ satisfies the properties given above Theorem \ref
{th-cs-non-percolation}.

For all $\alpha\ge0, \beta>0$ we define a random set $\Sigma(\alpha
,\beta)$ by
\[
\Sigma(\alpha,\beta) = \bigcup_{(c,r)\in\xi  \dvtx    r \in[\alpha
,\beta]} B(c,r).
\]
Notice that this set is empty if $\beta$ is strictly smaller than
$\alpha$.
If $x$ belongs to $\mathbb{R}^d$, we define an event $G(x,\alpha
,\beta)$ by
\[
G(x,\alpha,\beta) =
 \left\{
\matrix{\mbox{the connected component of } \Sigma(\alpha,\beta) \cup
B(x,\beta) \cr
\mbox{containing } x \mbox{ is not contained in }B(x,2\beta)}
\right\}.
\]
In other words, $G(x,\alpha,\beta)$ occurs if one can go from
$B(x,\beta)$ to the complement of $B(x,2\beta)$
using balls of the percolation process whose radii belong to $[\alpha
,\beta]$.
By stationarity of $\xi$, the probability of $G(x,\alpha,\beta)$
does not depend on $x$.
We denote it by $\pi(\alpha,\beta)$
\[
\pi(\alpha,\beta)=P(G(0,\alpha,\beta)).
\]

Similarly, for all $\beta>0$, we define an event $\widetilde{G}(\beta
)$ by
\[
\widetilde{G}(\beta)=
\left\{
\matrix{\mbox{the connected component of } \Sigma\cup B(0,\beta) \cr
\mbox{containing } 0 \mbox{ is not contained in } B(0,2\beta)}
\right\}.
\]
We denote its probability by $\widetilde{\pi}(\beta)$
\[
\widetilde{\pi}(\beta)=P(\widetilde{G}(\beta)).
\]

In order to state some relations between percolation and the various
events we have already introduced,
we shall need the following two events.
For all $\beta>0$ and $\rho>1$ we define $\widetilde{H}(\beta)$ and
$H(\rho,\beta)$ by
\[
\widetilde{H}(\beta) = \{\exists(c,r)\in\xi\dvtx  B(c,r)\cap B(0,2\beta
)\neq\varnothing\mbox{ and } r > \beta\}
\]
and
\[
H(\rho,\beta)=\{\exists(c,r) \in\xi\dvtx  c \in B(0,3\rho\beta)
\mbox{ and } r \in[\beta, \rho\beta]\}.
\]

We will give a strenghtened version of Theorem \ref
{th-cs-non-percolation} in which we relax the independence assumption.
To state this result, we shall need the following definition,
in which $\rho$ is strictly larger that $1$ and $\alpha, \beta$ are
as above, that is, $\alpha\ge0$ and $\beta>0$.
\begin{eqnarray*}
\hspace*{-2pt}&&I(\rho,\alpha,\beta) \\
\hspace*{-2pt}&&\qquad=
\sup_{x\in\mathbb{R}^d\setminus B(0,\rho\beta)}
\bigl[ P \bigl(G(0,\alpha,\beta) \cap G(x,\alpha,\beta) \bigr) -
P (G(0,\alpha,\beta) )P (G(x,\alpha,\beta) )  \bigr].
\end{eqnarray*}
Note that, under assumption B0 of Theorem \ref
{th-cs-non-percolation}, $I(\rho,\alpha,\beta) = 0$ for large enough~$\rho$
(see the beginning of the proof of Theorem \ref{th-cs-non-percolation}).
We also let
\[
I^+(\rho,\alpha,\beta) = \max (I(\rho,\alpha,\beta),0 ).
\]

\subsection{Proof of key inequalities}
\label{s-key}



Let us recall that $\mu$ is defined above the statement of Theorem
\ref{th-cs-non-percolation}
and that $M$ is defined in Section \ref{s-notations}.
The aim of this subsection is to prove the following result.

\begin{prop} \label{laprop}
Let $\rho\ge2$.
There exists a constant $\widetilde{D}>0$, that depends only on the
dimension $d$ and on $\rho$,
such that the following assertion holds for all $\alpha\ge0$ and all
$\beta>0$:
%
%
\begin{equation} \label{prop1}
\pi(\alpha,\rho\beta)\le
\widetilde{D}\pi(\alpha,\beta)^2
+\widetilde{D}\int_{[\beta,\rho\beta]} r^d\mu(dr)
+\widetilde{D}I^+(\rho,\alpha,\beta).
\end{equation}
Moreover, for all $\beta>0$, we have
%
%
\begin{equation} \label{prop3}
\pi(0,\beta)=\lim_{\alpha\to0} \pi(\alpha,\beta)
\end{equation}
and
%
%
\begin{equation} \label{prop2}
  P(M > 2\beta) \le\widetilde{\pi}(\beta) \le\pi(0,\beta
)+\widetilde{D} \int_{[\beta,+\infty[} r^d\mu(dr).
\end{equation}
\end{prop}

\begin{remark*} With (\ref{prop1}), we relate percolation probabilities at
different scales.
Our strategy is therefore related to multiscale strategies developed
for example in
\cite{Menshikovalmulti} and \cite{Menshikovalmultiunbounded}
(which use some stochastic domination properties) or in
\cite{FreirePopovVachkovskiapercolationmarriage}
(from which our approach is closer).
\end{remark*}

\subsubsection*{Ideas of the proof of the first part of Theorem
\protect\ref{th-cs-non-percolation} using Proposition \protect\ref{laprop}}
The aim is to prove that $P(M \ge\beta)$ tends to $0$ when $\beta$
tends to infinity.
By (\ref{prop2}) we get, under assumption B2, that it is
sufficient to prove that $\pi(0,\beta)$ tends to $0$.
By (\ref{prop3}) we get that it is sufficient to prove that $\pi
(\alpha,\beta)$ tends to $0$ uniformly in $\alpha$.
But by~(\ref{prop1}), $\pi(\alpha,\rho\beta)$ is bounded above by
$\widetilde{D}\pi(\alpha,\beta)^2$ up to error terms which satisfy
the following properties:
\begin{enumerate}
\item They are bounded above, by assumption B1;
\item They tend to $0$ when $\beta$ tends to infinity, by assumptions
B0 and B2.
\end{enumerate}
As $\pi(\alpha,\beta)=0$ for small enough $\beta$ (this is why the
parameter $\alpha$ has been introduced)
and as the bound given by assumption B1 on error terms is
small enough,
we first deduce that $\pi(\alpha,\beta)$ remains small for all
values of $\beta$ (see the first item of Lemma \ref{l-analyse}).
Then, as the error terms tend to $0$,
we get that $\pi(\alpha,\beta)$ tends to $0$ as~$\beta$ tends to
infinity (see the second item of Lemma \ref{l-analyse}).

The key lemma in the proof of Proposition \ref{laprop} is the
following one.
\begin{lemma} \label{rec}
Let $\rho\ge2$.
There exists a positive constant $D_1$ that depends only on the
dimension $d$ and on $\rho$ such that,
for all $\alpha\ge0$ and all $\beta>0$, the following holds:
\[
\pi(\alpha,\rho\beta) \le D_1\pi(\alpha,\beta)^2+D_1 I(\rho
,\alpha,\beta)+ P(H(\rho,\beta)).
\]
\end{lemma}
\begin{pf}
For all $r\ge0$ we denote by $S_r$ the Euclidean sphere
centered at the origin with radius $r$:
\[
S_r=\{x\in\mathbb{R}^d \dvtx  \|x\|=r\}.
\]
We fix $K$ and $L$, two finite subsets of $\mathbb{R}^d$ such that the
following properties hold:
\[
K \subset S_{\rho}\subset K+B(0,1) \quad\mbox{and}\quad L \subset S_{2\rho}
\subset L+B(0,1).
\]
We define $D_1$ as the product of the cardinalities of the sets $K$ and $L$.

Let $\alpha\ge0$ and $\beta>0$.
In this step, we prove the following inclusion:
%
%
\begin{equation} \label{hehe}
G(0,\alpha,\rho\beta) \setminus H(\rho,\beta) \subset
\biggl(\bigcup_{k\in K} G(\beta k,\alpha,\beta) \biggr) \cap
\biggl(\bigcup_{l\in L} G(\beta l,\alpha,\beta) \biggr).
\end{equation}

%
%
\begin{figure}

\includegraphics{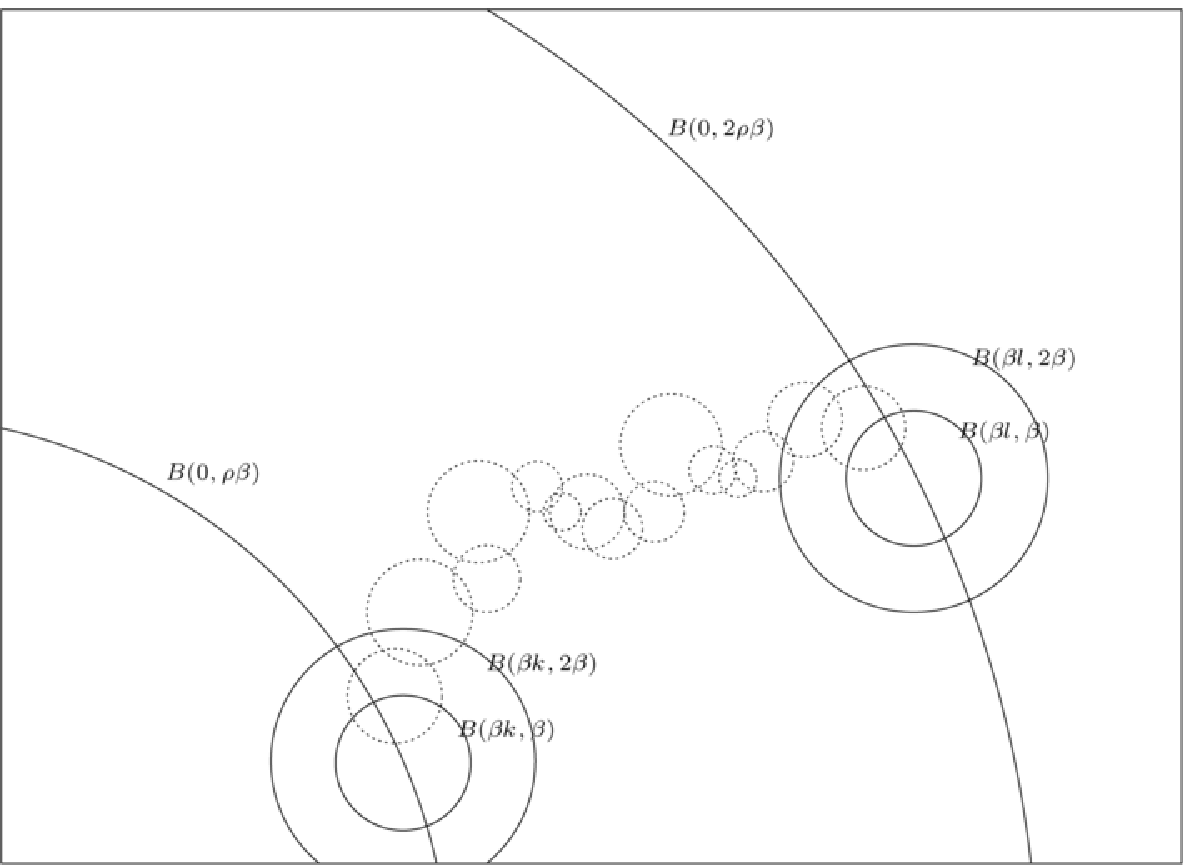}

\caption{Proof of (\protect\ref{hehe}).}
\label{zoli}
\end{figure}

We assume that the event $G(0,\alpha,\rho\beta)$ occurs but that the
event $H(\rho,\beta)$ does not occur.
As $G(0,\alpha,\rho\beta)$ occurs, one can go from $S_{\rho\beta}$
to $S_{2\rho\beta}$
using only balls of the percolation process whose radii belong to
$[\alpha,\rho\beta]$.
One can furthermore assume that the center of each such ball belongs to
$B(0,\rho3\beta)$.

We refer to Figure \ref{zoli} where the doted circles stand for some of the
previous balls.
One of these balls touches $S_{\rho\beta}$.
This ball then touches $B(\beta k, \beta)$ for some $k\in K$.
We then see that one can go from $B(\beta k, \beta)$ to the complement
of $B(\beta k, 2\beta)$
using only balls whose radii belong to $[\alpha,\rho\beta]$ and
whose centers belong to $B(0,3\rho\beta)$.

But, as $H(\rho,\beta)$ does not occur, the radius of each such ball
$B(c,r)$ is less than $\beta$.
Therefore, $G(\beta k, \alpha, \beta)$ occurs.
We have proved that the event $\bigcup_{k\in K} G(\beta k, \alpha,
\beta)$ occurs.
We can prove in a similar way that the event $\bigcup_{l\in L} G(\beta
l,\alpha,\beta)$ occurs.
Therefore the inclusion (\ref{hehe}) is proved.

We then get from (\ref{hehe})
\[
\pi(\alpha,\rho\beta)
\le P(H(\rho,\beta))+\sum_{k\in K, l\in L} P\bigl(G(\beta k,\alpha,\beta
)\cap G(\beta l,\alpha,\beta)\bigr).
\]
For all $k\in K$ and all $l\in L$, we have $\|\beta k -\beta l\| \ge
\beta\rho$.
By stationarity and by definition of $I(\rho,\alpha,\beta)$ and of
$D_1$, we then get
\[
\pi(\alpha,\rho\beta)
\le P(H(\rho,\beta))+D_1  \bigl(\pi(\alpha,\beta)^2+I(\rho,\alpha
,\beta) \bigr).
\]
This ends the proof.
\end{pf}
\begin{lemma} \label{isildur}
For all $\beta>0$, the following holds:
\[
\pi(0,\beta)=\lim_{\alpha\to0} \pi(\alpha,\beta).
\]
\end{lemma}
\begin{pf}
Let $\beta>0$.
As $\alpha\mapsto\Sigma(\alpha,\beta)$ is nonincreasing, $\alpha
\mapsto G(0,\alpha,\beta)$ is nonincreasing.
Consequently,
\[
\lim_{\alpha\to0} \pi(\alpha,\beta)
=P \biggl(\bigcup_{\alpha>0} G(0,\alpha,\beta) \biggr).
\]
Therefore, it is sufficient to prove the following equality:
\[
\bigcup_{\alpha>0} G(0,\alpha,\beta)=G(0,0,\beta).
\]
If the event $G(0,0,\beta)$ occurs, then one can go from $B(0,\beta)$
to the complement of $B(0,2\beta)$
using balls of the percolation process whose radii belongs to $]0,\beta]$.
By a compactness argument, we get the existence of a real $\alpha>0$
such that
one can go from $B(0,\beta)$ to the complement of $B(0,2\beta)$
using balls of the percolation process whose radii belongs to $[\alpha
,\beta]$.
In other words, $G(0,\alpha,\beta)$ occurs.
This proves one of the required inclusions.
The other inclusion is straightforward.
\end{pf}
\begin{lemma} \label{isidore}
For all $\beta>0$, the following inclusion holds:
\[
\{M > 2\beta\} \subset\widetilde{G}(\beta) \subset G(0,0,\beta)
\cup\widetilde{H}(\beta).
\]
\end{lemma}
\begin{pf}
%
%
Let $\beta>0$.
If $G(0,0,\beta)$ does not occur, then one cannot go from $B(0,\beta
)$ to the complement of $B(0,2\beta)$
using balls\vspace*{1pt} of the percolation process whose radii belongs to $]0,\beta]$.
If moreover $\widetilde{H}(\beta)$ does not occur,
then balls of the percolation process whose radii do not belong to
$]0,\beta]$
will not help to connect $B(0,\beta)$ to the complement of $B(0,2\beta)$.
Therefore $\widetilde{G}(\beta)$ does not occur.
This proves one inclusion.
The other one is straightforward.
\end{pf}
\begin{lemma} \label{ph}
There exists a positive constant $D_2$, that depends only on the
dimension $d$, such that
for all $\beta>0$,
the following inequality holds:
\[
P(\widetilde{H}(\beta)) \le D_2 \int_{[\beta,+\infty[} r^d \mu(dr).
\]
\end{lemma}
\begin{pf}
We have
\[
\widetilde{H}(\beta)=\{\xi\cap V(\beta) \neq\varnothing\},
\]
where
\[
V(\beta)=\{(c,r) \in\mathbb{R}^d\times\,]0,+\infty[ \dvtx  B(c,r)\cap
B(0,2\beta
)\neq\varnothing\mbox{ and } r>\beta\}.
\]
We therefore have
\begin{eqnarray*}
P(\widetilde{H}(\beta))
& = & P\bigl(\xi\cap V(\beta) \neq\varnothing\bigr) \\
& \le& E\bigl(\operatorname{card}\bigl(\xi\cap V(\beta)\bigr)\bigr) \\
& = & \int_{\mathbb{R}^d} dc \int_{]0,+\infty[} \mu(dr) 1_{V(\beta)}(c,r).
\end{eqnarray*}
As
\begin{eqnarray*}
V(\beta)
& = & \{(c,r) \in\mathbb{R}^d\times\,]0,+\infty[ \dvtx  \|c\|<r+2\beta
\mbox{ and
} r>\beta\}
\end{eqnarray*}
we get
\begin{eqnarray*}
P(\widetilde{H}(\beta))
& \le& \int_{]\beta,+\infty[} |B(0,r+2\beta)| \mu(dr) \\
& \le& \int_{]\beta,+\infty[} |B(0,3r)| \mu(dr).
\end{eqnarray*}
The inequality stated in the lemma is therefore fulfilled with $D_2=|B(0,3)|$.
\end{pf}
\begin{lemma} \label{pht}
Let $\rho\ge2$.
There exists a positive constant $D_3$, that depends only on the
dimension $d$ and on $\rho$, such that
for all $\beta>0$,
the following inequality holds:
\[
P(H(\rho,\beta)) \le D_3 \int_{[\beta,\rho\beta]} r^d \mu(dr).
\]
\end{lemma}
\begin{pf}
We have
\begin{eqnarray*}
P(H(\rho,\beta))
& \le& E \bigl(\operatorname{card} \bigl(\{(c,r)\in\xi\dvtx  c \in
B(0,3\rho\beta)
\mbox{ and } r \in[\beta, \rho\beta]\} \bigr) \bigr) \\
& = & |B(0,3\rho\beta)| \mu([\beta, \rho\beta]) \\
& = & |B(0,3\rho)|\beta^d \mu([\beta, \rho\beta]).
\end{eqnarray*}
The inequality stated in the lemma is therefore fulfilled with
$D_3=|B(0,3\rho)|$.
\end{pf}

\begin{pf*}{Proof of Proposition \protect\ref{laprop}}
This is a consequence of the previous lemmas.
Let us denote by $D_1$, $D_2$ and $D_3$ the constants given by Lemmas
\ref{rec}, \ref{ph} and \ref{pht}.
Set $\widetilde{D}=\max(1,D_1,D_2,D_3)$.
By Lemma \ref{rec} we have, for all $\alpha\ge0$ and $\beta>0$
\[
\pi(\alpha,\rho\beta) \le\widetilde{D}\pi(\alpha,\beta
)^2+\widetilde{D} I^+(\rho,\alpha,\beta)+ P(H(\rho,\beta)).
\]
By Lemma \ref{pht} we then get (\ref{prop1}).
By Lemma \ref{isidore} we have, for all $\beta>0$
\[
P(M > 2\beta) \le\widetilde{\pi}(\beta) \le\pi(0,\beta) +
P(\widetilde{H}(\beta)).
\]
By Lemma \ref{ph} we then get (\ref{prop2}).
Finally, (\ref{prop3}) is given by Lemma \ref{isildur}.
\end{pf*}

\subsection[Proof of Theorem 1.3]{Proof of Theorem \protect\ref{th-cs-non-percolation}}
\label{s.existence}

We first state three theorems which, together,
give a strenghtened version of Theorem \ref{th-cs-non-percolation}.
Notice that the conclusion of each of the first two theorems is one of
the assumptions of the following one.
\begin{theorem} \label{th-cs-non-percolation-I}
For all $\rho\ge2$ and $D>0$, consider the following hypothesis
$\mathcal{H}(\rho,D)$.

There exist sequences ${(\alpha_n)}_{n\in\mathbb{N}}$ and ${(\beta
_n)}_{n\in\mathbb{N}}$ of real numbers such that the following
conditions hold:
\begin{enumerate}
\item For all $n\in\mathbb{N}$, $\alpha_n \ge0$, and $\alpha_n$
tends to
$0$ when $n$ tends to infinity.
\item For all $n\in\mathbb{N}$, $\beta_n > 0$, and ${(\beta
_n)}_{n\in\mathbb{N}}$
is bounded.
\item For all $n\in\mathbb{N}$ and all $\beta\ge\beta_n$,
$I^+(\rho
,\alpha_n,\beta) \le D$ and $\beta^d \mu([\beta,+\infty[) \le D$.
\item For all $n\in\mathbb{N}$ and all $\beta\in[\beta_n,\rho
\beta_n]$,
$\pi(\alpha_n,\beta) \le D$.
\end{enumerate}
Then, for all $\rho\ge2$ and $D'>0$,
there exists $D>0$, that depends only on $d$, $\rho$ and~$D'$, such that
$\mathcal{H}(\rho,D)$ implies that the probability $\pi(0,\beta)$
is smaller than $D'$ for large enough $\beta$.
\end{theorem}
\begin{remarks*}
\begin{enumerate}
\item
Note that, in the first assumption, one allows the sequence to be
constant equal to $0$.
\item If $\beta$ belongs to $]0,\alpha_n[$ then $\Sigma(\alpha
_n,\beta)$ is empty and therefore $G(0,\alpha_n,\beta)$ cannot occur.
The probability $\pi(\alpha_n,\beta)$ then equals $0$.
Therefore, the fourth assumption of $\mathcal{H}(\rho,D)$ is always
satisfied when $\beta_n$ is strictly smaller than $\alpha_n\rho^{-1}$.
This is the reason why we introduced the parameter $\alpha$ in the
definition of $G$ and $\pi$.
\end{enumerate}
\end{remarks*}
\begin{theorem} \label{th-cs-non-percolation-II}
For all $\rho\ge2$ and $D'>0$, consider the following hypothesis
$\mathcal{H'}(\rho,D')$:
\begin{enumerate}
\item The probability $\pi(0,\beta)$ is smaller than $D'$ for large
enough $\beta$.
\item$I^+(\rho,0,\beta)$ tends to $0$ as $\beta$ tends to infinity.
\item The integral $\int_{[1,+\infty[} \beta^d \mu(d\beta)$ is finite.
\end{enumerate}
Let $\rho\ge2$.
There exists $D'>0$, that depends only on $d$ and $\rho$, such that
$\mathcal{H'}(\rho,D')$ implies
that the probability $\widetilde{\pi}(\beta)$ tends to $0$ as $\beta
$ tends to infinity
and, therefore, implies that there is almost surely no percolation.
\end{theorem}
\begin{theorem} \label{th-cs-non-percolation-III}
Let $\rho\ge2$ and $s>0$.
Assume the following:
\begin{enumerate}
\item The probability $\widetilde{\pi}(\beta)$ tends to $0$ as
$\beta$ tend to infinity.
\item$\int_{[1,+\infty[} \beta^{s-1} I^+(\rho,0,\beta) \,d\beta
<\infty$.
\item$\int_{[1,+\infty[} \beta^{d+s} \mu(d\beta)<\infty$.
\end{enumerate}
Then, the integral
\[
\int_0^{+\infty} \beta^{s-1} \widetilde{\pi}(\beta)\,d\beta
\]
is finite.
Therefore, the moment $E(M^s)$ is finite.
\end{theorem}

The proof of the previous theorems relies on Proposition \ref{laprop}
and on the following elementary lemma.
There are three items in the lemma.
Each of them corresponds to one of the previous theorems.
\begin{lemma} \label{l-analyse}
Let $f$ and $g$ be two measurable functions from $]0,+\infty[$ to
$[0,+\infty[$.
Let $\rho>1$.
We assume that, for all $\beta>0$, the following inequality holds:
%
%
\begin{equation} \label{hou}
f(\rho\beta)\le f(\beta)^2+g(\beta).
\end{equation}
Then:
\begin{enumerate}
\item Let $\varepsilon\in\,]0,1]$.
If there exists $\beta_0>0$ such that $f(\beta)\le\varepsilon/2$ for
all $\beta\in[\beta_0,\rho\beta_0]$
and $g(\beta)\le\varepsilon/4$ for all $\beta\ge\beta_0$ then, for
all $\beta\ge\beta_0$, we have $f(\beta)\le\varepsilon/2$.
\item If, for all large enough $\beta>0$, the inequality $f(\beta)\le
1/2$ holds and if
$g(\beta)$ converges to $0$ as $\beta$ tends to infinity then,
$f(\beta)$ converges to $0$ as $\beta$ tends to infinity.
\item Let $s>-1$ be a real number.
If $f$ is bounded, if $f(\beta)$ converges to $0$ as $\beta$ tends to
infinity and
if the integral $\int_1^{+\infty} \beta^s g(\beta) \,d\beta$ is
finite then,
the integral $\int_0^{+\infty} \beta^s f(\beta) \,d\beta$ is finite.
\end{enumerate}
\end{lemma}
\begin{pf}
\textit{Proof of item $1$}.\quad If $\beta>0$ is such that $f(\beta
)\le\varepsilon/2$ and $g(\beta)\le\varepsilon/4$, then
\[
f(\rho\beta)\le\varepsilon^2/4+\varepsilon/4 \le\varepsilon/2.
\]
The result follows.

\textit{Proof of item $2$}.\quad By (\ref{hou}) we get
\[
\limsup_{\beta\to\infty} f(\beta) \le \Bigl[\limsup_{\beta\to
\infty} f(\beta) \Bigr]^2+\limsup_{\beta\to\infty} g(\beta).
\]
By assumption,
\[
\limsup_{\beta\to\infty} f(\beta) \le1/2 \quad\mbox{and}\quad \limsup
_{\beta\to\infty} g(\beta)=0.
\]
As $f$ is nonnegative, we get that $f(\beta)$ converges to $0$ as
$\beta$ tends to infinity.

%
%

\textit{Proof of item $3$}.\quad
Let $s>-1$.
As $f$ tends to $0$, there exists a real $A\ge\rho$ such that
%
%
\begin{equation} \label{xien}
\forall\beta\ge A\rho^{-1} \dvtx  f(\beta)\le\rho^{-s-1}/2.
\end{equation}
For all real $r\ge A$, we get, by (\ref{hou}) and (\ref{xien})
\begin{eqnarray*}
&& \int_A^r f(\beta)\beta^s\,d\beta\\
&&\qquad \le \int_A^r f(\beta\rho^{-1})^2\beta^s \,d\beta+ \int_A^r
g(\beta\rho^{-1})\beta^s\,d\beta\\
&&\qquad \le \rho^{s+1} \int_{A\rho^{-1}}^{r\rho^{-1}} f(\beta)^2\beta
^s \,d\beta
+ \rho^{s+1}\int_{A\rho^{-1}}^{+\infty} g(\beta)\beta^s\,d\beta\\
&&\qquad \le 1/2 \int_{A\rho^{-1}}^{r\rho^{-1}} f(\beta)\beta^s \,d\beta
+ \rho^{s+1}\int_{A\rho^{-1}}^{+\infty} g(\beta)\beta^s\,d\beta\\
&&\qquad \le 1/2 \int_{A}^r f(\beta)\beta^s \,d\beta
+ 1/2 \int_{A\rho^{-1}}^A f(\beta)\beta^s \,d\beta
+ \rho^{s+1}\int_{A\rho^{-1}}^{+\infty} g(\beta)\beta^s\,d\beta.
\end{eqnarray*}
As $f$ is bounded, the integral $\int_A^r f(\beta)\beta^s \,d\beta$
is finite.
We therefore get
\[
\int_A^r f(\beta)\beta^s\,d\beta\le\int_{A\rho^{-1}}^A f(\beta
)\beta^s \,d\beta
+ 2\rho^{s+1}\int_{A\rho^{-1}}^{+\infty} g(\beta)\beta^s\,d\beta
\]
and then
\[
\int_A^{+\infty} f(\beta)\beta^s\,d\beta\le\int_{A\rho^{-1}}^A
f(\beta)\beta^s \,d\beta
+ 2\rho^{s+1}\int_1^{+\infty} g(\beta)\beta^s\,d\beta.
\]
As $f$ is bounded, the lemma follows.
\end{pf}
\begin{pf*}{Proofs of Theorems \protect\ref{th-cs-non-percolation-I},
\protect\ref
{th-cs-non-percolation-II} and \protect\ref{th-cs-non-percolation-III}}
Let $\widetilde{D}$ be the positive constant given by Proposition \ref
{laprop}.
For all $\alpha\ge0$ we define a function $f_{\alpha}\dvtx ]0,+\infty
[\,\to[0,+\infty[$ by
\[
f_{\alpha}(\beta)=\widetilde{D}\pi(\alpha,\beta)
\]
and a function $g_{\alpha}\dvtx ]0,+\infty[\,\to[0,+\infty[$ by
\[
g_{\alpha}(\beta)= \widetilde{D}^2 I^+(\rho,\alpha,\beta
)+\widetilde{D}^2 \int_{[\beta,\rho\beta]} r^d\mu(dr).
\]
By (\ref{prop1}) we get, for all $\alpha\ge0$ and all $\beta>0$
%
%
\begin{equation} \label{appliprop1}
f_{\alpha}(\rho\beta)\le f_{\alpha}(\beta)^2+g_{\alpha}(\beta).
\end{equation}

\textit{Proof of Theorem \ref{th-cs-non-percolation-I}.}\quad
Let
\[
\varepsilon=\min \bigl(\tfrac12,2D'\widetilde{D} \bigr)>0
\]
and
\[
D=\min \biggl(\frac{\varepsilon}{8\rho^2\widetilde{D}^2},\frac
{\varepsilon}{2\widetilde{D}} \biggr)>0.
\]
Let us prove that $D$ satisfies the required properties of Theorem \ref
{th-cs-non-percolation-I}.
Let ${(\alpha_n)}_n$ and ${(\beta_n)}_n$ be as in the statement of
the theorem.
Let $\beta_*$ be the supremum of the bounded sequence ${(\beta_n)}_n$.

Let $n\in\mathbb{N}$.
By the third assumption of hypothesis $\mathcal{H}(\rho,D)$ we get,
for all \mbox{$\beta\ge\beta_n$},
\begin{eqnarray*}
g_{\alpha_n}(\beta)
& \le& \widetilde{D}^2\rho^d\beta^d\mu([\beta,+\infty
[)+\widetilde{D}^2I^+(\rho,\alpha_n,\beta) \\
& \le& 2\widetilde{D}^2\rho^d D \\
& \le& \varepsilon/4.
\end{eqnarray*}
By the fourth assumption of hypothesis $\mathcal{H}(\rho,D)$ we get,
for all $\beta\in[\beta_n,\rho\beta_n]$,
\[
f_{\alpha_n}(\beta) \le\widetilde{D}D \le\varepsilon/2.
\]
By the first item of Lemma \ref{l-analyse}, we then get the inequality
$f_{\alpha_n}(\beta)\le\varepsilon/2$
for all $\beta\ge\beta_n$.
Therefore, for all $\beta\ge\beta_*$, we have
\[
\pi(\alpha_n,\beta) \le D'.
\]
The theorem follows thanks to Lemma \ref{isildur}.

\textit{Proof of Theorem
\ref{th-cs-non-percolation-II}.}\quad
Let
\[
D'=\frac1{2\widetilde{D}}>0.
\]
Let us check that $D'$ satisfies the required properties of Theorem
\ref{th-cs-non-percolation-II}.
By the first assumption of the theorem, we know that the inequality
$\pi(0,\beta)\le D'$ holds for large enough $\beta$.
Therefore, we have $f_0(\beta)\le1/2$ for large enough $\beta$.
By the second and the third assumptions, we get that $g_0(\beta)$
converges to $0$ as $\beta$ tends to infinity.
By the second item of Lemma \ref{l-analyse}, we then get that
$f_0(\beta)$ also converges to $0$.
Therefore, $\pi(0,\beta)$ converges to $0$.
The theorem follows thanks to the third assumption and to (\ref{prop2}).

\textit{Proof of theorem
\ref{th-cs-non-percolation-III}.}\quad
For all $\beta>0$, we have $\pi(0,\beta)\le\widetilde{\pi}(\beta)$.
By the first assumption of the theorem, we then have the convergence of
$\pi(0,\beta)$ to $0$.
Therefore, $f_0(\beta)$ converges to $0$.
Let us notice the following:
%
%
\begin{eqnarray}\label{youpi}
\int_1^{+\infty} d\beta\,\beta^{s-1} \int_{[\beta,+\infty[} \mu
(dr) r^d
& = & \int_{[1,\infty[} \mu(dr) r^d \int_1^r d\beta\,\beta^{s-1}
\nonumber\\
& \le& \int_{[1,\infty[} \mu(dr) s^{-1}r^{d+s} \\
& < & \infty\nonumber
\end{eqnarray}
by the third assumption.
Using also the second assumption, we then get that the integral $\int
_1^{+\infty} \beta^{s-1}g_0(\beta)\,d\beta$ is finite.
By the third item of Lemma \ref{l-analyse},
we then get that the integral $\int_0^{+\infty} \beta^{s-1}
f_0(\beta)\,d\beta$ is finite.
The integral
%
%
\begin{equation} \label{gym}
\int_0^{+\infty} \beta^{s-1} \pi(0,\beta)\,d\beta
\end{equation}
is therefore also finite.
But by (\ref{prop2}) we have, for all $\beta>0$
\[
\widetilde{\pi}(\beta) \le\pi(0,\beta)+\widetilde{D}\int
_{[\beta,+\infty[}r^d\mu(dr).
\]
By (\ref{youpi}) and (\ref{gym}), we thus get that the integral
$\int_1^{+\infty} \beta^{s-1}\widetilde{\pi}(\beta)\,d\beta$
and then the integral $\int_0^{+\infty} \beta^{s-1}\widetilde{\pi
}(\beta)\,d\beta$ is finite.
The theorem follows by the first inequality of~(\ref{prop2}).
\end{pf*}

\begin{pf*}{Proof of Theorem \protect\ref{th-cs-non-percolation}}
Let $\rho=\max(4C,2)$.
Let $D'$ be the constant given by Theorem \ref{th-cs-non-percolation-II}.
Let $D$ be the constant given by Theorem \ref{th-cs-non-percolation-I}.
Let us check that $D$ satisfies the required properties.

Let $\alpha\ge0$ and $\beta>0$.
Let us notice that, for all $x\in\mathbb{R}^d$, the event $G(x,\alpha
,\beta
)$ only depends on
$\xi\cap B(x,3\beta)\times\,]0,\beta]$.
Therefore, the event $G(x,\alpha,\beta)$ only depends on $\xi\cap
B(x,3\beta)\times\,]0,3\beta]$.
By assumption B0 we then get that $G(0,\alpha,\beta)$ and
$G(x,\alpha,\beta)$ are independent
as soon as $\|x\| \ge3\beta C$.
By definition of $\rho$, we thus get
%
%
\begin{equation} \label{Inul}
I(\rho,\alpha,\beta)=0.
\end{equation}

Let $n$ be a positive integer.
We let $\alpha_n=n^{-1}$ and $\beta_n=\alpha_n(2\rho)^{-1}$.
For all $\beta\in[\beta_n,\rho\beta_n]$, $\beta$ belongs to
$]0,\alpha_n[$.
Therefore the set $\Sigma(\alpha_n,\beta)$ is empty and consequently
the event $G(0,\alpha_n,\beta)$ does not occur.
As a consequence, the fourth assumption of hypothesis $\mathcal{H(\rho
,D)}$ holds.

The third assumption of hypothesis $\mathcal{H}(\rho,D)$ holds
because of (\ref{Inul}) and assumption B1 of Theorem \ref
{th-cs-non-percolation}.
By Theorem \ref{th-cs-non-percolation-I}, we then get that $\pi
(0,\beta)$ is smaller than~$D'$ for large enough $\beta$.
In other words, the first assumption of Theorem \ref
{th-cs-non-percolation-II} holds.
The second assumption of this theorem holds because of (\ref{Inul}).
The third one holds because of assumption B2 of Theorem \ref
{th-cs-non-percolation}.
We then get that $S$ is almost surely bounded and that the first
assumption of Theorem \ref{th-cs-non-percolation-III} holds.
By (\ref{Inul}), the second assumption of Theorem \ref
{th-cs-non-percolation-III} holds.
If assumption B3 holds, we then get, by Theorem \ref
{th-cs-non-percolation-III}, that $E(M^s)$ is finite.
\end{pf*}

\section[Proofs of Theorems 1.1 and 1.2 (models induced by Poisson
processes)]{Proofs of Theorems \protect\ref{th-poisson-non-percolation}
and \protect\ref{th-poisson-moment-fini} (models induced by Poisson
processes)}
\label{s-preuve-poisson}

We work with the objects defined in Section \ref{s-poisson}.
In particular, $\xi$ is a Poisson point process on $\mathbb
{R}^d\times
\,]0,+\infty[$ and we have
\[
\xi=\{(c,r(c)), c\in\chi\},
\]
where $\chi$ denotes the projection of $\xi$ on $\mathbb{R}^d$.

The following elementary lemma is stated and proven in \cite
{Gbooleanmodel} for a probability measure $\mu$.
The proof is the same for a locally finite measure.
\begin{lemma} \label{l-perco-grand}
Let $\mu$ be a locally finite measure on $]0,+\infty[$.
If $\int_{[1,+\infty[}\beta^d\times\break\mu(d\beta)$ is infinite then, for
all $\lambda>0$,
we have $P_{\lambda,\mu}$-almost surely $\Sigma=\mathbb{R}^d$.
If $s>0$ is such that $\int_{[1,+\infty[}\beta^{d+s}\mu(d\beta)$
is infinite then, for all $\lambda>0$, $E_{\lambda,\mu}(M^s)$ is infinite.
\end{lemma}
\begin{pf}
Let $\mu$ be a locally finite measure on $]0,+\infty[$ and
$\lambda>0$.

We first prove that, for all $r>0$, the following inequality holds:
%
%
\begin{eqnarray}\label{pilou}
&&P_{\lambda,\mu} \bigl(\exists c\in\chi\dvtx  B(0,r) \subset B(c,r(c)) \bigr)
\nonumber\\[-8pt]\\[-8pt]
&&\qquad\ge1 - \exp \biggl(-\lambda2^{-d}|B(0,1)| \int_{[2r,+\infty[} \beta
^d \mu(d\beta)  \biggr).\nonumber
\end{eqnarray}
Let $r>0$. We have
\[
P_{\lambda,\mu}\bigl(\exists c\in\chi\dvtx  B(0,r) \subset B(c,r(c))\bigr) = P(\xi
\cap A\neq\varnothing),
\]
where
\[
A=\{(c,\beta) \in\xi\dvtx  \beta\ge\|c\|+r\}.
\]
Therefore
\begin{eqnarray*}
P_{\lambda,\mu}\bigl(\exists c\in\chi\dvtx  B(0,r) \subset B(c,r(c))\bigr)
& = & 1 - \exp \biggl(-\lambda\int_{\mathbb{R}^d} \mu([\|c\|
+r,+\infty[)
\,dc  \biggr) \\
& = & 1 - \exp \biggl(-\lambda\int_{[r,+\infty[} |B(0,\beta-r)| \mu
(d\beta) \biggr) \\
& \ge& 1 - \exp \biggl(-\lambda\int_{[2r,+\infty[} [B(0,\beta-r)|
\mu(d\beta) \biggr) \\
& \ge& 1 - \exp \biggl(-\lambda\int_{[2r,+\infty[} |B(0,\beta
/2)|\mu(d\beta) \biggr).
\end{eqnarray*}
The relation (\ref{pilou}) is proved.

If $\int_{[1,+\infty[} \beta^d \mu(d\beta)$ is infinite
then, by (\ref{pilou}), we get, for all $r>0$
\[
P_{\lambda,\mu} \bigl(\exists c\in\chi\dvtx  B(0,r) \subset B(c,r(c))
\bigr) = 1.
\]
Therefore, almost surely, we have $\Sigma=\mathbb{R}^d$.

Let $s>0$.
We assume now that $\int_{[1,+\infty[}\beta^{d+s}\mu(d\beta)$ is infinite.
If $\int_{[1,+\infty[} \beta^d\times\break \mu(d\beta)$ is infinite,
the desired result is a consequence of what we have proved in the
previous step.
We assume henceforth that $\int_{[1,+\infty[} \beta^d \mu(d\beta)$
is finite.
Let $C$ be defined by
\[
C=\lambda2^{-d}|B(0,1)| \int_{[1,+\infty[} \beta^d \mu(d\beta).
\]
This constant is finite.
By (\ref{pilou}) we get, for all $r>1/2$, the following inequality:
\begin{eqnarray*}
&&P_{\lambda,\mu} \bigl(\exists c\in\chi\dvtx  B(0,r) \subset B(c,r(c)) \bigr)
\\
&&\qquad\ge C^{-1}\bigl(1-\exp(-C)\bigr) \lambda2^{-d}|B(0,1)| \int_{[2r,+\infty[}
\beta^d \mu(d\beta)
\end{eqnarray*}
and then
%
%
\begin{equation} \label{clementine}
P(M \ge r) \ge C^{-1}\bigl(1-\exp(-C)\bigr) \lambda2^{-d}|B(0,1)| \int
_{[2r,+\infty[} \beta^d \mu(d\beta).
\end{equation}
As $\int_{[1,+\infty[}\beta^{d+s}\mu(d\beta)$ is infinite, the integral
\[
\int_{1/2}^{+\infty}  \biggl(r^{s-1} \int_{[2r,+\infty[} \beta^d
\mu(d\beta) \biggr)\,dr
\]
is infinite.
Therefore, by (\ref{clementine}), the integral
$\int_0^{+\infty} r^{s-1} P_{\lambda,\mu}(M\ge r) \,dr$ is infinite.
The moment $E_{\lambda,\mu}(M^s)$ is then infinite.
\end{pf}
\begin{lemma} \label{l-Abis}
Let $\mu$ be a locally finite measure on $]0,+\infty[$.
Let $\rho>1$.
We have
\[
\sup_{r>0} r^d \mu([r,\rho r]) \le\sup_{r>0} r^d \mu([r,+\infty[)
\le\frac{1}{1-\rho^{-d}} \sup_{r>0} r^d\mu([r,\rho r]).
\]
\end{lemma}
\begin{pf}
The first inequality is straightforward. Let us prove the other one.
Let $r>0$.
We have
\begin{eqnarray*}
r^d \mu([r,+\infty[)
& = & \sum_{n\ge0} \rho^{-nd} (r\rho^n)^d\mu([r\rho^n, r\rho
^{n+1}[) \\
& \le& \sum_{n\ge0} \rho^{-nd} \sup_{s>0} s^d\mu([s, s\rho[)\\
& \le& \frac{1}{1-\rho^{-d}} \sup_{s>0} s^d\mu([s, s\rho[).
\end{eqnarray*}
The lemma follows.
\end{pf}
\begin{lemma} \label{l-perco-petit}
Assume $d \ge2$.
Let $\mu$ be a locally finite measure on $]0,+\infty[$.
If
\[
\sup_{r>0} r^d\mu([r,+\infty[)
\]
is infinite, then for all $\lambda>0$,
we have $P_{\lambda,\mu}(\mbox{percolation})>0$.
\end{lemma}
\begin{pf}
Let $\mu$ be a locally finite measure on $]0,+\infty[$ and
$\lambda>0$.
Let $\lambda_c>0$ be the critical value for the classical Boolean
model when all radii equal $1$
(see, e.g., \cite{MeesterRoylivre} or Section 12.10 in \cite
{Grimmettpercolation}).
In other words, when $\mu=\delta_1$, $S$ is almost surely bounded
when $\lambda<\lambda_c$
and $S$ is unbounded with positive probability when $\lambda>\lambda_c$.

Let $\rho=2$.
By assumption and by Lemma \ref{l-Abis}, there exists $r_0>0$ such that
\[
\lambda r_0^d\mu([r_0,r_0\rho]) > \lambda_c.
\]
We define a new Poisson point process as follows:
\[
\widetilde{\xi}=\{(c,r_0) \dvtx  c\in\chi\mbox{ such that } r(c)\in
[r_0,r_0\rho]\}.
\]
The\vspace*{1pt} intensity measure of this point process is the product of the
measure $\lambda\mu([r_0,r_0\rho])|\cdot|$ by
the probability measure $\delta_{r_0}$.
Let $\widetilde{\Sigma}$ be associated with $\widetilde{\xi}$ as in
Section \ref{s-notations}.
Let us notice that $\widetilde{\Sigma}$ is a subset of $\Sigma$.
It is therefore sufficient to prove that $\widetilde{\Sigma}$ is in
the supercritical phase.
The random set $r_0^{-1}\widetilde{\Sigma}$ is associated with the
following Poisson point process
\[
r_0^{-1}\widetilde{\xi}=\{(cr_0^{-1},1) \dvtx  c\in\chi\mbox{ such that
} r(c)\in[r_0,r_0\rho]\}
\]
whose intensity measure is the product of $r_0^d\lambda\mu
([r_0,r_0\rho])|\cdot|$ by the probability measure $\delta_1$.
By our choice of $r_0$ and by definition of $\lambda_c$ we get that
$r_0^{-1}\widetilde{\Sigma}$,
and therefore $\widetilde{\Sigma}$, is in the supercritical phase.
This ends the proof.
\end{pf}

\begin{pf*}{Proofs of Theorems \protect\ref
{th-poisson-non-percolation} and \protect\ref
{th-poisson-moment-fini}}

\textit{Proof of sufficient conditions.}\quad
Let $C=2$.
Let $D>0$ be the constant given by Theorem \ref{th-cs-non-percolation}.
Assumption B0
of Theorem \ref{th-cs-non-percolation} is satisfied because of
independence properties of
Poisson point processes.
Since, under $P_{\lambda,\mu}$,
the intensity measure of $\xi$ is the product of the Lebesgue measure
and of the measure $\lambda\mu$,
the required results follow from Theorem
\ref{th-cs-non-percolation}.

\textit{Proof of necessary conditions.}\quad
This is a consequence of Lemmas~\ref{l-perco-petit} and \ref
{l-perco-grand}.
\end{pf*}

\section[Proofs of Lemma 1.4 and Theorem 1.8 (multiscale percolation)]{Proofs of
Lemma \protect\ref{l-fractal} and Theorem \protect\ref{th-fractal} (multiscale
percolation)}
\label{s-preuve-fractal}

\mbox{}

\begin{pf*}{Proof of Lemma \protect\ref{l-fractal}}
Let us first notice that, for each $n\ge0$, $a^{-n}\xi_n$ is
a~Poisson point process whose intensity measure
is the product of $\lambda|\cdot|$ by the measure $\nu_n$ defined by
$\nu_n(B)=a^{nd}\nu(a^n B)$.
Let us recall that the measure $\mu$ was defined in (\ref{f-def-mu}) by
\[
\mu(B)=\sum_{n \ge0} a^{nd} \nu(a^n B).
\]
We then have $\mu=\sum_n \nu_n$.

It remains to check that the measure $\mu$ is locally finite.
Let $k\in\mathbb{Z}$.
It is sufficient to prove that $\mu([a^k,a^{k+1}[)$ is finite.
We have
\begin{eqnarray*}
\mu([a^k,a^{k+1}[)
& = & \sum_{n \ge0} a^{nd} \nu([a^{k+n},a^{k+n+1}[) \\
& \le& \int_{[a^k,+\infty[} x^d a^{-kd} \nu(dx).
\end{eqnarray*}
As $\int_{]0,+\infty[} x^d \nu(dx)$ is finite, the result follows.
\end{pf*}

\begin{pf*}{Proof of Theorem \protect\ref{th-fractal}}
Let $\mu$ be the measure defined by (\ref{f-def-mu}).
Thanks to Theorem \ref{th-poisson-non-percolation} is it sufficient to
check the following:
\begin{enumerate}
\item Condition A1 holds.
\item Condition A2 holds if and only if $\int
_{[1,+\infty[} \beta^d\ln(\beta)\nu(d\beta)$ is finite.
\end{enumerate}
Let us notice that, for all $f\dvtx ]0,+\infty[\,\to\mathbb{R}$ measurable and
nonnegative, we have
\[
\int_{]0,+\infty[}f(\beta)\mu(d\beta)=\sum_{n\ge0} a^{nd}\int
_{]0,+\infty[}f(a^{-n}\beta)
\nu(d\beta).
\]

Let us check the first item.
Let $r>0$.
We have
\begin{eqnarray*}
\int_{[r,ra]} \beta^d \mu(d\beta)
& = & \sum_{n\ge0} a^{nd} \int_{]0,+\infty[} 1_{[r,ra]}(\beta
a^{-n}) (\beta a^{-n})^d \nu(d\beta) \\
& = & \int_{]0,+\infty[} \sum_{n\ge0} 1_{[r,ra]}(\beta a^{-n})
\beta^d \nu(d\beta) \\
& \le& \int_{]0,+\infty[} 2\beta^d \nu(d\beta).
\end{eqnarray*}
The first item then follows from Lemma \ref{l-Abis} by (\ref
{hyp-moment-fractal}).

Let us check the second item.
As above, we get
\begin{eqnarray*}
\int_{[1,+\infty[} \beta^d \mu(d\beta)
& = & \sum_{n\ge0} a^{nd} \int_{]0,+\infty[} 1_{[1,+\infty[}(\beta
a^{-n}) (\beta a^{-n})^d \nu(d\beta) \\
& = & \int_{]0,+\infty[} \sum_{n\ge0} 1_{[1,+\infty[}(\beta
a^{-n}) \beta^d \nu(d\beta) \\
& = & \int_{[1,+\infty[}  \bigl( \lfloor\ln(\beta)\ln(a)^{-1}
 \rfloor+1 \bigr) \beta^d \nu(d\beta).
\end{eqnarray*}
The second item follows.
This concludes the proof.
\end{pf*}

\section[Proofs of Lemma 1.11 and Theorem 1.10 (marriage)]{Proofs of
Lemma \protect\ref{l-marriage-domination} and
Theorem \protect\ref{th-marriage-non-perco-moments} (marriage)}
\label{s-proof-marriage}

Let us recall the definition of $\xi$.
We assume that $\chi$ is a Poisson point process on $\mathbb{R}^d$ whose
intensity measure is the Lebesgue measure.
For all $a\in\chi$ we define $R(a,\chi)$ by
\[
R(a,\chi)=\inf \bigl\{r \ge0 \dvtx  \alpha\operatorname{card}\bigl(\chi\cap
\overline
{B}(a,2r)\bigr) \le|B(a,r)| \bigr\}.
\]
[We let $R(a,\chi)=\infty$ if there is no such $r$.]
Using some elementary properties of the map defined by $r\mapsto\alpha
\operatorname{card}(\chi\cap\overline{B}(a,2r))-|B(a,r)|$, we get
that $R(a,\chi)$ is always positive and that
\[
R(a,\chi)=\min \bigl\{r \ge0\dvtx
|B(0,r)|\in\alpha\mathbb{N}\quad\mbox{and}\quad
\alpha\operatorname{card}\bigl(\chi\cap\overline{B}(a,2r)\bigr) =
|B(a,r)| \bigr\}.
\]
(With the same convention as before if there is no such $r$.)
Among other things, this remark enables us to easily solve some
measurability issues.
We define a point process $\xi$ on $\mathbb{R}^d\times\,]0,+\infty]$ by
\[
\xi=\{(a,2R(a,\chi)), a\in\chi\}.
\]
Let us notice that the law of $\xi$ is invariant under the action of
the translations of~$\mathbb{R}^d$
and that the intensity measure of $\xi$ is locally finite.
The intensity measure is therefore the product of the Lebesgue measure
on $\mathbb{R}^d$
by a locally finite measure on $]0,+\infty]$.
We denote this measure on $]0,+\infty]$ by $\mu$.
Let us notice that $\mu$ is a probability measure.
\begin{lemma} \label{l-marriage-queue-mu}
There exists an absolute constant $K>0$ and a function $F\dvtx
]0,2^{-d}[\,\to\,
]0,+\infty[$ that depends only on the dimension $d$ such that:
\begin{enumerate}
\item $\lim_{\alpha\to0} F(\alpha)=+\infty$.
\item For all $\alpha\in\,]0,2^{-d}[$ and all $r>0$, we have: $\mu
(]r,+\infty]) \le K\exp(-F(\alpha)r^d).$
\end{enumerate}
\end{lemma}
\begin{pf}
Assume $\alpha\in\,]0,2^{-d}[$.
Let $r>0$.
By definition of $\mu$ and $\xi$ we have
\begin{eqnarray*}
\mu(]r,+\infty])
& = & E \bigl(\operatorname{card}(\xi\cap[0,1]^d\times\,]r,+\infty
]) \bigr) \\
& = & E \biggl(\sum_{a\in\chi\cap[0,1]^d} 1_{2R(a,\chi) > r}
\biggr) \\
& = & E \biggl(\sum_{a\in\chi\cap[0,1]^d} 1_{2R(0,\chi-a) >
r} \biggr).
\end{eqnarray*}
As the Palm measure of the Poisson point process $\chi$ is the law of
$\chi\cup\{0\}$ (see, e.g., \cite{Moller}), we get
\[
\mu(]r,+\infty]) = P\bigl(2R(0,\chi\cup\{0\}) > r\bigr).
\]
By definition of $R(0,\chi\cup\{0\})$, we then get
\begin{eqnarray*}
\mu(]r,+\infty])
& \le& P \bigl(\alpha\operatorname{card} \bigl((\chi\cup\{0\}) \cap
\overline
{B}(0,r) \bigr) > |B(0,r/2)| \bigr) \\
& = & P \bigl(\alpha\bigl(N(r)+1\bigr) > \omega_d r^d2^{-d} \bigr) \\
& = & P \bigl(N(r) > \alpha^{-1}\omega_d r^d2^{-d}-1 \bigr),
\end{eqnarray*}
where
$N(r)=\operatorname{card}(\chi\cap\overline{B}(0,r))$ and $\omega
_d=|B(0,1)|$.

If
%
%
\begin{equation} \label{torduoups}
1<\alpha^{-1}\omega_d r^d2^{-d} \bigl(1-\sqrt{\alpha2^d}\bigr)
\end{equation}
we have
\[
\mu(]r,+\infty]) \le P \bigl(N(r) > \omega_d r^d\sqrt{\alpha
2^d}^{-1} \bigr).
\]
As $N(r)$ is a Poisson random variable with mean $w_dr^d$ we then get,
using Chernoff's bound
\[
\mu(]r,+\infty]) \le\exp\bigl(-w_d r^d g\bigl(\sqrt{\alpha2^d}\bigr)\bigr),
\]
where $g\dvtx ]0,1[\,\to\mathbb{R}$ is defined by
\[
g(x)=\bigl(x-1-\ln(x)\bigr)/x.
\]
The previous inequality holds as soon as (\ref{torduoups}) holds.
It therefore holds as soon as $\omega_d r^d > \alpha2^d (1-\sqrt
{\alpha2^d}) ^{-1}$.

Now, if $\omega_d r^d \le\alpha2^d (1-\sqrt{\alpha2^d}) ^{-1}$
then, as $g$ is nonnegative
\[
\omega_d r^d g\bigl(\sqrt{\alpha2^d}\bigr) \le
h\bigl(\sqrt{\alpha2^s}\bigr),
\]
where $h\dvtx ]0,1[\,\to\mathbb{R}$ is defined by
\[
h(x)=x^2(1-x)^{-1}g(x).
\]
As $\mu(]r,+\infty])$ is at most $1$ we have, in this case
%
%
\begin{equation} \label{ecole}
\mu(]r,+\infty]) \le\exp\bigl(h\bigl(\sqrt{\alpha2^d}\bigr)\bigr)\exp\bigl(-w_d r^d
g\bigl(\sqrt{\alpha2^d}\bigr)\bigr).
\end{equation}

As $h$ is nonnegative, we finally get that (\ref{ecole}) holds for
all $r>0$.
As $h$ is bounded and as $\lim_{x\to0} g(x)=+\infty$, the lemma follows.
\end{pf}

We assume henceforth that $\alpha$ is strictly smaller than $2^{-d}$.
By the previous lemma, we can therefore consider that $\xi$ is a point
process on $\mathbb{R}^d\times\,]0,+\infty[$
and that $\mu$ is a probability measure on $]0,+\infty[$.
We are therefore in the same framework as in Section \ref{s-general}.
We associate with $\xi$ a random set $\Sigma$ and a random variable~$M$.

\begin{pf*}{Proof of Lemma \protect\ref{l-marriage-domination}}
We work on a full event on which there exists an a.e. unique stable
allocation and denote by $\psi$ one of those allocations.
Let $a\in\chi$.
Let us recall that $R(a,\chi)$ is finite.
To simplify notation, we write $R$ instead of $R(a,\chi)$.
To prove the lemma, it suffices to check that $\psi^{-1}(a)$ is a
subset of $\overline{B}(a,R)$.
We have
\[
\alpha\operatorname{card}\bigl(\chi\cap\overline{B}(a,2 R)\bigr) = |B(a,R)|.
\]
Let $\varepsilon>0$ be such that there is no point of $\chi$ in the
shell $\overline{B}(a,2R+2\varepsilon) \setminus\overline{B}(a,2R)$.
We then have
\[
\alpha\operatorname{card}\bigl(\chi\cap\overline{B}(a,2R + 2\varepsilon
)\bigr) <
|B(a,R+\varepsilon)|.
\]
Therefore
\[
 \bigl| \psi^{-1} \bigl( \chi\cap\overline{B}(a,2R + 2\varepsilon
) \bigr)  \bigr| < |B(a,R+\varepsilon)|.
\]
As a consequence, there exists $x$ in $B(a,R+\varepsilon)$ such that
$\psi(x)$ belongs to $\chi\cup\{\infty\}$
and does not belong to $\overline{B}(a,2R+2\varepsilon)$.
If $\psi(x)\in\chi$, we have
\[
\|x-\psi(x)\|>R+\varepsilon\quad\mbox{and}\quad \|x-a\| \le R+\varepsilon.
\]
In particular, $x$ desires $a$.
Otherwise, that is, if $\psi(x)=\infty$, then $x$ also desires $a$.
As $\psi$ is stable, we therefore get that $a$ does not covet $x$.
As a consequence, $\psi^{-1}(a)$ is contained in $\overline{B}(a,\|
x-a\|)$ and therefore in $\overline{B}(a,R+\varepsilon)$.
As this result holds for arbitrary small $\varepsilon>0$, we get that
$\psi^{-1}(a)$ is contained in $\overline{B}(a,R)$.
The lemma follows.
\end{pf*}

\begin{pf*}{Proof of Theorem \protect\ref{th-marriage-non-perco-moments}}
Thanks to Lemma \ref{l-marriage-domination},
it suffices to check that~$\xi$ satisfies the assumptions of Theorem
\ref{th-cs-non-percolation}.
\begin{enumerate}[B2.]
\item[B0.] We show that the assumption is fulfiled with $C=7$.
Let $r>0$.
For all $a\in\chi$ we let
\[
\widetilde{R}(a,\chi)=\inf\bigl\{s \in[0,r] \dvtx  \alpha\operatorname
{card}\bigl(\chi\cap
\overline{B}(a,2s)\bigr) \le|B(a,s)|\bigr\}.
\]
[We let $\widetilde{R}(a,\chi)=r$ if there exists no such $s$.]
Let us notice that, for all $a\in\chi$, we have $\widetilde
{R}(a,\chi)=R(a,\chi)$
as soon as $R(a,\chi)<r$ or $\widetilde{R}(a,\chi)<r$.
Therefore, for all $x\in\mathbb{R}^d$,
\[
\xi\cap\mathbb{R}^d \times[0,r[\, = \widetilde{\xi} \cap\mathbb
{R}^d\times[0,r[,
\]
where $\widetilde{\xi}$ is defined by
\[
\widetilde{\xi}=\{(a,2\widetilde{R}(a,\chi)), a\in\chi\}.
\]
As a consequence, we see that $\xi\cap B(x,r) \times[0,r[$ only
depends on $\chi\cap B(x,3r)$.
By the independence property of Poisson point processes, we then get
that, if $x$ belongs to $\mathbb{R}^d\setminus B(0,6r)$,
the point processes $\xi\cap B(0,r) \times[0,r[$ and $\xi\cap B(x,r)
\times[0,r[$ are independent.
The required result follows.
\item[B1.] By Lemma \ref{l-marriage-queue-mu}, we have
\begin{eqnarray*}
\sup_{r>0} r^d \mu([r,+\infty[)
& \le& \sup_{r>0} r^d \mu(]r/2,+\infty[) \\
& \le& \sup_{r>0} r^d K\exp(-F(\alpha)r^d2^{-d}) \\
& = & K2^dF(\alpha)^{-1} \sup_{x>0} x\exp(-x).
\end{eqnarray*}
As $F(\alpha)$ tends to infinity when $\alpha$ tends to $0$,
assumption B1 is fulfiled for small enough $\alpha$.
\item[B2 and B3.]
By Lemma \ref{l-marriage-queue-mu}, we get that $\int_{]0,+\infty[}
r^{d+s} \mu(dr)$ is finite for all $s \ge0$.
\end{enumerate}
When $\alpha$ is small enough, we can thus use Theorem \ref
{th-cs-non-percolation}.
We get that $E(M^s)$ is finite for all $s>0$.
By Lemma \ref{l-marriage-domination} we then get that $E(D^s)$ is
finite for all $s>0$.
\end{pf*}

\section*{Acknowledgment}
I would like to thank the referees for helpful remarks.


%
\printaddresses


\begin{thebibliography}{18}

\bibitem{FreirePopovVachkovskiapercolationmarriage}
\begin{barticle}[msn]
\bauthor{\bsnm{Freire},~\bfnm{M.~V.}\binits{M.~V.}},
  \bauthor{\bsnm{Popov},~\bfnm{S.}\binits{S.}} \AND
  \bauthor{\bsnm{Vachkovskaia},~\bfnm{M.}\binits{M.}}
(\byear{2007}).
\btitle{Percolation for the stable marriage of {P}oisson and {L}ebesgue}.
\bjournal{Stochastic Process. Appl.}
\bvolume{117}
\bpages{514--525}.
\bmrnumber{MR2305384}
\end{barticle}
\endbibitem

\bibitem{GaleShapleymarriage}
\begin{barticle}[msn]
\bauthor{\bsnm{Gale},~\bfnm{D.}\binits{D.}} \AND
  \bauthor{\bsnm{Shapley},~\bfnm{L.~S.}\binits{L.~S.}}
(\byear{1962}).
\btitle{College admissions and the stability of marriage}.
\bjournal{Amer. Math. Monthly}
\bvolume{69}
\bpages{9--15}.
\bmrnumber{MR1531503}
\end{barticle}
\endbibitem

\bibitem{Gbooleanmodel}
\begin{barticle}[msn]
\bauthor{\bsnm{Gou{\'e}r{\'e}},~\bfnm{Jean-Baptiste}\binits{J.-B.}}
(\byear{2008}).
\btitle{Subcritical regimes in the {P}oisson {B}oolean model of continuum
  percolation}.
\bjournal{Ann. Probab.}
\bvolume{36}
\bpages{1209--1220}.
\bmrnumber{MR2435847}
\end{barticle}
\endbibitem

\bibitem{Grimmettpercolation}
\begin{bbook}[msn]
\bauthor{\bsnm{Grimmett},~\bfnm{Geoffrey}\binits{G.}}
(\byear{1999}).
\btitle{Percolation},
\bedition{2nd} ed.
\bseries{Grundlehren der Mathematischen Wissenschaften [Fundamental Principles
  of Mathematical Sciences]}
\bvolume{321}.
\bpublisher{Springer}, \baddress{Berlin}.
\bmrnumber{MR1707339}
\end{bbook}
\endbibitem

\bibitem{Hallcontinuumpercolation}
\begin{barticle}[msn]
\bauthor{\bsnm{Hall},~\bfnm{Peter}\binits{P.}}
(\byear{1985}).
\btitle{On continuum percolation}.
\bjournal{Ann. Probab.}
\bvolume{13}
\bpages{1250--1266}.
\bmrnumber{MR806222}
\end{barticle}
\endbibitem

\bibitem{Halllivre}
\begin{bbook}[vtex]
\bauthor{\bsnm{Hall},~\bfnm{Peter}\binits{P.}}
(\byear{1988}).
\btitle{Introduction to the Theory of Coverage Processes}.
\bpublisher{Wiley}, \baddress{New York}.
\bmrnumber{MR973404}
\end{bbook}
\endbibitem

\bibitem{HoffmanHolroydPeresmarriage}
\begin{barticle}[msn]
\bauthor{\bsnm{Hoffman},~\bfnm{Christopher}\binits{C.}},
  \bauthor{\bsnm{Holroyd},~\bfnm{Alexander~E.}\binits{A.~E.}} \AND
  \bauthor{\bsnm{Peres},~\bfnm{Yuval}\binits{Y.}}
(\byear{2006}).
\btitle{A stable marriage of {P}oisson and {L}ebesgue}.
\bjournal{Ann. Probab.}
\bvolume{34}
\bpages{1241--1272}.
\bmrnumber{MR2257646}
\end{barticle}
\endbibitem

\bibitem{HoffmanHolroydPeresmarriagetail}
\begin{bmisc}[unstr]
\bauthor{\bsnm{Hoffman},~\bfnm{Christopher}\binits{C.}},
  \bauthor{\bsnm{Holroyd},~\bfnm{Alexander E.}\binits{A. E.}} \AND
  \bauthor{\bsnm{Peres},~\bfnm{Yuval}\binits{Y.}}
  (\byear{2008}).
\bhowpublished{Tail bounds for the stable marriage of Poisson and Lebesgue.
\textit{Canad. J. Math.} To appear. Available at}
\href{http://www.arxiv.org/abs/math/0507324}{arXiv:math/0507324}.
\end{bmisc}
\endbibitem


\bibitem{Kallenbergrandommeasures}
\begin{bbook}[msn]
\bauthor{\bsnm{Kallenberg},~\bfnm{Olav}\binits{O.}}
(\byear{1986}).
\btitle{Random Measures}, \bedition{4th} ed.
\bpublisher{Akademie-Verlag}, \baddress{Berlin}.
\bmrnumber{MR854102}
\end{bbook}
\endbibitem

\bibitem{MeesterRoylivre}
\begin{bbook}[msn]
\bauthor{\bsnm{Meester},~\bfnm{Ronald}\binits{R.}} \AND
  \bauthor{\bsnm{Roy},~\bfnm{Rahul}\binits{R.}}
(\byear{1996}).
\btitle{Continuum Percolation}.
\bseries{Cambridge Tracts in Mathematics}
\bvolume{119}.
\bpublisher{Cambridge Univ. Press}, \baddress{Cambridge}.
\bmrnumber{MR1409145}
\end{bbook}
\endbibitem

\bibitem{MeesterRoySarkar}
\begin{barticle}[msn]
\bauthor{\bsnm{Meester},~\bfnm{Ronald}\binits{R.}},
  \bauthor{\bsnm{Roy},~\bfnm{Rahul}\binits{R.}} \AND
  \bauthor{\bsnm{Sarkar},~\bfnm{Anish}\binits{A.}}
(\byear{1994}).
\btitle{Nonuniversality and continuity of the critical covered volume fraction
  in continuum percolation}.
\bjournal{J. Statist. Phys.}
\bvolume{75}
\bpages{123--134}.
\bmrnumber{MR1273055}
\end{barticle}
\endbibitem

\bibitem{Menshikovalmulti}
\begin{barticle}[msn]
\bauthor{\bsnm{Menshikov},~\bfnm{M.~V.}\binits{M.~V.}},
  \bauthor{\bsnm{Popov},~\bfnm{S.~Yu.}\binits{S.~Y.}} \AND
  \bauthor{\bsnm{Vachkovskaia},~\bfnm{M.}\binits{M.}}
(\byear{2001}).
\btitle{On the connectivity properties of the complementary set in fractal
  percolation models}.
\bjournal{Probab. Theory Related Fields}
\bvolume{119}
\bpages{176--186}.
\bmrnumber{MR1818245}
\end{barticle}
\endbibitem

\bibitem{Menshikovalmultiunbounded}
\begin{barticle}[msn]
\bauthor{\bsnm{Menshikov},~\bfnm{M.~V.}\binits{M.~V.}},
  \bauthor{\bsnm{Popov},~\bfnm{S.~Yu.}\binits{S.~Y.}} \AND
  \bauthor{\bsnm{Vachkovskaia},~\bfnm{M.}\binits{M.}}
(\byear{2003}).
\btitle{On a multiscale continuous percolation model with unbounded defects}.
\bjournal{Bull. Braz. Math. Soc. (N.S.)}
\bvolume{34}
\bpages{417--435}.
\bmrnumber{MR2045167}
\end{barticle}
\endbibitem

\bibitem{MenshikovSidorenkocoincidence}
\begin{barticle}[vtex]
\bauthor{\bsnm{Menshikov},~\bfnm{M.~V.}\binits{M.~V.}} \AND
  \bauthor{\bsnm{Sidorenko},~\bfnm{A.~F.}\binits{A.~F.}}
(\byear{1987}).
\btitle{Coincidence of critical points in {P}oisson percolation models}.
\bjournal{Teor. Veroyatnost. i Primenen.}
\bvolume{32}
\bpages{603--606}.
\bmrnumber{MR914958}
\end{barticle}
\endbibitem

\bibitem{Moller}
\begin{bbook}[msn]
\bauthor{\bsnm{M{\o}ller},~\bfnm{Jesper}\binits{J.}}
(\byear{1994}).
\btitle{Lectures on Random {V}orono\u\i\ Tessellations}.
\bseries{Lecture Notes in Statistics}
\bvolume{87}.
\bpublisher{Springer}, \baddress{New York}.
\bmrnumber{MR1295245}
\end{bbook}
\endbibitem

\bibitem{Neveupp}
\begin{bincollection}[vtex]
\bauthor{\bsnm{Neveu},~\bfnm{J.}\binits{J.}}
(\byear{1977}).
\btitle{Processus ponctuels}.
In \bbooktitle{\'{E}cole D'\'{E}t\'e de {P}robabilit\'es de {S}aint-{F}lour,
  {VI}---1976}.
\bseries{Lecture Notes in Math.}
\bvolume{598}
\bpages{249--445}.
\bpublisher{Springer}, \baddress{Berlin}.
\bmrnumber{MR0474493}
\end{bincollection}
\endbibitem

\bibitem{ZuevI}
\begin{barticle}[msn]
\bauthor{\bsnm{Zuev},~\bfnm{S.~A.}\binits{S.~A.}} \AND
  \bauthor{\bsnm{Sidorenko},~\bfnm{A.~F.}\binits{A.~F.}}
(\byear{1985}).
\btitle{Continuous models of percolation theory. {I}}.
\bjournal{Teoret. Mat. Fiz.}
\bvolume{62}
\bpages{76--86}.
\bmrnumber{MR782099}
\end{barticle}
\endbibitem

\bibitem{ZuevII}
\begin{barticle}[msn]
\bauthor{\bsnm{Zuev},~\bfnm{S.~A.}\binits{S.~A.}} \AND
  \bauthor{\bsnm{Sidorenko},~\bfnm{A.~F.}\binits{A.~F.}}
(\byear{1985}).
\btitle{Continuous models of percolation theory. {II}}.
\bjournal{Teoret. Mat. Fiz.}
\bvolume{62}
\bpages{253--262}.
\bmrnumber{MR783056}
\end{barticle}
\endbibitem

\end{thebibliography}
\end{document}